\documentstyle[12pt]{article}

\begin{document}

\date{}

\title{\Large\bf Reformulating the Map Color Theorem}

\author{Louis H. Kauffman \\ Department of Mathematics, Statistics and Computer
Science \\ University of Illinois at Chicago \\ 851 South Morgan Street\\
Chicago, IL, 60607-7045}

\maketitle

\thispagestyle{empty}

\subsection*{\centering Abstract}

{\em This paper discusses reformulations of the problem of coloring plane maps with
 four colors.  We include discussion of the Eliahou-Kryuchkov conjecture, the Penrose
 formula, the vector cross product formulation and
 the reformulations in terms of formations and factorizations due to G. Spencer-Brown.}

\section{Introduction}

In this paper we give a concise introduction to the work of G. Spencer-Brown
\cite{SB} on the Four-Color Theorem and some of the consequences of this work
in relation to other reformulations of the four color problem.  
This work involves a rewriting of the
coloring problem in terms of two-colored systems of Jordan curves in the plane. 
These systems, called {\em formations}, are in one-to-one correspondence with
cubic plane graphs that are colored with three edge colors so that three distinct
colors are incident to each vertex of the graph. It has long been known that the
four color problem can be reformulated in terms of coloring such cubic graphs.
\vspace{3mm}

We first concentrate on proving two key results.  The first is a Parity Lemma
due to G. Spencer-Brown \cite{SB}. This lemma is also implied by work of Tutte \cite{T} via
translation from edge colorings to formations.
The second result, depending on the Parity
Lemma is a proof that a certain principle of irreducibility for formations is
{\em equivalent} to the Four-Color  Theorem.  Spencer-Brown takes this principle
of irreducibility (here called the {\em Primality Principle}) to be axiomatic and
hence obtains a proof of the Four-Color Theorem that is based upon it. He also
gives proofs of the Primality Principle (See Theorem 17 \cite{SB}, p. 168-170)
that depend upon a subtle notion of inverse distinction. This work of Spencer-Brown deserves careful
consideration. \vspace{3mm}

The present paper is an expansion of \cite{LK}.  In that paper we also  prove the parity lemma 
and discuss the primality principle. However, the discussion of factorizability of formations
is imprecise in \cite{LK} and I have taken the opportunity of this paper to rectify that fault.
I hope that this paper attains the desired clarity in regard to parity and primality. In the 
author's opinion these concepts are central to understanding the nature of the four color theorem, 
and it is worth a second try at explication.
\bigbreak

There are seven sections to the present paper. In the second section we give the
basics about cubic maps and formations. In the third section we prove the Parity
Lemma. In the fourth section we give the equivalence of the Four-Color Theorem and
the Primality Principle. In the fifth section we discuss an algorithm, the Parity Pass, discovered
by Spencer-Brown.  The Parity Pass is an algorithm designed to color a map that has been colored
except for a five sided region. The language of the algorithm is in terms of formations. It is an
extraordinarily powerful algorithm and may in itself constitute a solution to the four color problem.
It is worth conjecturing that this is so. In the sixth section we discuss an application of 
formations to the workings of a chromatic counting formula due to Roger Penrose. In section seven we 
apply ideas from formations to the Eliahou-Kryuchkov conjecture, showing that it can be 
reformulated in terms of coloring and re-coloring trees, and in terms of the vector cross product
reformulation of the four color theorem.  
\bigbreak

\noindent {\bf Acknowledgement.}    It gives the author pleasure to thank
James Flagg and Karanbir Sarkaria for helpful conversations in the course of 
constructing this paper.
\bigbreak

\section{\bf Cubic Graphs and Formations}

A {\em graph} consists in a vertex set $V$ and an edge set $E$ such that every
edge has two vertices associated with it (they may be identical).  If a vertex is
in the set of vertices associated with an edge, we say that this vertex {\em
belongs} to that edge. If two vertices form the vertex set for a given edge we
say that that edge {\em connects} the two vertices (again the two may be
identical).  A {\em loop} in a graph is an edge whose vertex set has cardinality
one. In a {\em multi-graph} it is allowed that there may be a multiplicity of
edges connecting a given pair of vertices. All graphs in this paper are
multi-graphs, and we shall therefore not use the prefix "multi" from here on.
\vspace{3mm}

A {\em cubic graph} is a graph in which every vertex either belongs to three
distinct edges, or there are two edges at the vertex with one of them a loop. A
{\em coloring} (proper coloring) of a cubic graph $G$ is an assignment of the
labels $r$ (red), $b$ (blue), and $p$ (purple) to the edges of the graph so that
three distinct labels occur at every vertex of the graph. This means that there
are three distinct edges belonging to each vertex and that it is possible to
label the graph so that three distinct colors occur at each vertex. Note that a
graph with a loop is not colorable. \vspace{3mm}

The simplest uncolorable cubic graph is illustrated in Figure 1.  For obvious
reasons, we refer to this graph as the {\em dumbell}.  Note that the dumbell is
planar. \vspace{3mm}

{\tt    \setlength{\unitlength}{0.92pt} \begin{picture}(143,61) \thinlines   
\put(94,30){\circle*{8}} \put(50,31){\circle*{8}} \put(49,31){\line(1,0){45}}
\put(113,31){\circle{40}} \put(30,30){\circle{40}} \end{picture}}

\begin{center} {\bf Figure 1 --- The Simplest Uncolorable Cubic Graph}
\end{center} \vspace{3mm}

An edge in a connected plane graph is said to be an {\em isthmus} if the deletion of that
edge results in a disconnected graph. It is easy to see that a connected plane
cubic graph without isthmus is loop-free. \vspace{3mm}

Heawood reformulated the four-color conjecture (which we will henceforth refer to
as the {\em Map Theorem}) for plane maps to a corresponding statement about the
colorability of plane cubic graphs. In this form the theorem reads \vspace{2mm}

\noindent {\bf Map Theorem for Cubic Graphs.}  A connected plane cubic graph without
isthmus is properly edge-colorable with three colors. \vspace{3mm}

We now introduce a diagrammatic representation for the coloring of a cubic graph.
Let $G$ be a cubic graph and let $C(G)$ be a coloring of $G.$ Using the colors
$r$, $b$ and $p$ we will write purple as a formal product of red and blue: $$p =
rb.$$

\noindent One can follow single colored paths on the coloring $C(G)$ in the
colors red and blue. Each red or blue path will eventually return to its starting
point, creating a circuit in that color. The red circuits are disjoint from one
another, and the blue circuits are disjoint from one another. Red circuits
and blue circuits may meet along edges in $G$ that are colored purple
($p=rb$).   In the case of a plane graph $G$, a meeting of two circuits may take
the form of one circuit crossing the other in the plane, or one circuit may share
an edge with another circuit, and then leave on the same side of that other circuit.
We call these two
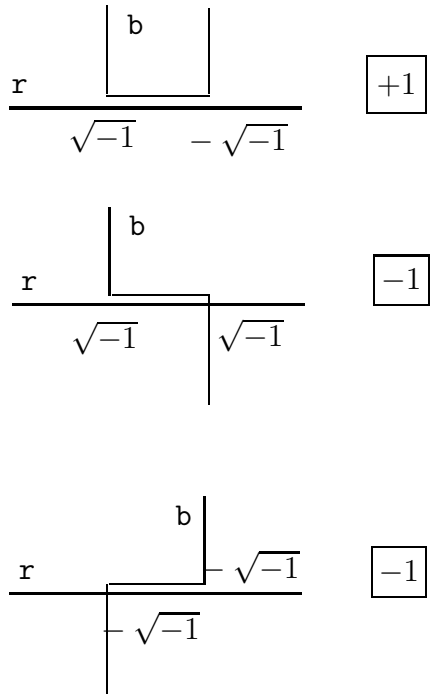
planar configurations a {\em cross} and a {\em bounce} respectively. \vspace{3mm}

{\tt    \setlength{\unitlength}{0.92pt} \begin{picture}(362,214) \thinlines   
\put(45,10){\makebox(67,18){formation}} \put(211,101){\makebox(67,18){graph}}
\put(196,95){\line(0,-1){30}} \put(244,62){\line(0,-1){29}}
\put(287,38){\makebox(51,17){bounce}} \put(204,12){\makebox(63,19){cross}}
\put(321,63){\line(0,1){29}} \put(296,63){\line(1,0){25}}
\put(296,92){\line(0,-1){29}} \put(227,64){\line(0,-1){29}}
\put(265,63){\line(-1,0){21}} \put(265,93){\line(0,-1){30}}
\put(196,64){\line(1,0){30}} \thicklines   \put(182,60){\line(1,0){170}}
\thinlines    \put(302,150){\makebox(26,23){G(F)}}
\put(130,58){\makebox(20,23){F}} \put(231,151){\makebox(17,18){b}}
\put(276,153){\makebox(17,17){p}} \put(181,152){\makebox(20,19){r}}
\put(77,37){\line(1,0){45}} \put(122,101){\line(0,-1){63}}
\put(76,101){\line(1,0){46}} \put(76,101){\line(0,-1){64}} \thicklines  
\put(35,31){\framebox(94,75){}} \put(247,199){\line(0,-1){74}}
\put(200,124){\framebox(94,75){}} \thinlines    \put(247,198){\circle*{12}}
\put(248,126){\circle*{12}} \put(10,122){\framebox(147,78){}}
\put(127,126){\makebox(17,18){p=br}} \put(108,148){\makebox(18,17){r}}
\put(107,176){\makebox(17,15){b}} \put(61,125){\makebox(59,18){purple}}
\put(22,131){\line(1,0){36}} \thicklines   \put(21,136){\line(1,0){36}}
\put(65,149){\makebox(32,16){red}} \put(60,174){\makebox(41,19){blue}}
\put(20,159){\line(1,0){34}} \thinlines    \put(21,183){\line(1,0){34}}
\end{picture}}

\begin{center} {\bf Figure 2 --- Coloring and Formation} \end{center}
\vspace{3mm}

\noindent {\bf Definition.} A {\em formation} \cite{SB} is a finite collection of 
simple closed curves in the plane, with each curve colored either red or blue such that the
red curves are disjoint from one another, the blue curves are disjoint from one
another and red and blue curves can meet in a finite number of segments (as
described above for the circuits in a coloring of a cubic graph). \vspace{2mm}

\noindent Associated with any formation $F$ there is a well-defined cubic graph
$G(F)$,  obtained by identifying the shared segments in the formation as edges in
the graph, and the endpoints of these segments as vertices. The remaining
(unshared) segments of each simple closed curve constitute the remaining edges of
$G(F).$  A formation $F$ is said to be a formation for a cubic graph $G$ if $G =
G(F).$ We also say that $F$ {\em formates} $G.$ \vspace{2mm}

\noindent A {\em plane formation} is a formation such that each simple closed
curve in the formation is a Jordan curve in the plane. For a plane formation,
each shared segment between two curves of different colors is either a bounce or
a crossing (see above), that condition being determined by the embedding of the
formation in the plane. \vspace{3mm}

Since the notion of a formation is abstracted from the circuit decomposition of a
colored cubic graph, we have the proposition: \vspace{3mm}

\noindent {\bf Proposition.}  Let $G$ be a cubic graph and $Col(G)$ be the set of
colorings of $G$.  Then $Col(G)$ is in one-to-one correspondence with the set of
formations for $G.$ 
\bigbreak

\noindent In particular, the Map Theorem is equivalent to the
\bigbreak

\noindent {\bf Formation Theorem.}  Every connected plane cubic graph without
isthmus has a formation. \vspace{3mm}

This equivalent version of the Map Theorem is due to G. Spencer-Brown \cite{SB}. 
The advantage of the Formation Theorem is that, just as one can enumerate graphs,
one can enumerate formations.  In particular, plane formations are generated by
drawing systems of Jordan curves in the plane that share segments according to
the rules explained above. This gives a new way to view the evidence for the Map
Theorem, since one can enumerate formations and observe that all the plane cubic
graphs are occurring in the course of the enumeration!   See Figures 2 and 3 for
illustrations of the relationship of formation with coloring. \vspace{3mm}

\noindent {\bf Remark.}  In the figures the reader will note that graphs are
depicted with horizontal and vertical edges. This means that some edges have
corners. These corners, artifacts of this form of representation, are not
vertices of the graph.  In depicting formations, we have endeavored to keep the
shared segments slightly separated for clarity in the diagram. These separated
segments are amalgamated in the graph that corresponds to the formation.
\vspace{3mm}

{\tt    \setlength{\unitlength}{0.92pt} \begin{picture}(385,181) \thinlines   
\put(256,10){\makebox(33,29){G}} \put(35,14){\makebox(31,27){F}} \thicklines  
\put(302,159){\line(1,0){48}} \thinlines    \put(325,84){\circle*{8}}
\put(301,54){\circle*{8}} \put(303,130){\circle*{8}} \put(325,73){\circle*{8}}
\put(284,55){\circle*{8}} \put(280,129){\circle*{8}} \thicklines  
\put(373,169){\line(0,-1){157}} \put(300,130){\line(1,0){24}}
\put(350,84){\line(-1,0){26}} \put(350,158){\line(0,-1){72}}
\put(280,130){\line(1,0){22}} \put(324,129){\line(0,-1){44}}
\put(302,158){\line(0,-1){28}} \put(324,72){\line(0,1){16}}
\put(284,55){\line(1,0){15}} \put(282,55){\line(-1,0){48}}
\put(283,74){\line(0,-1){19}} \put(323,74){\line(-1,0){39}}
\put(324,57){\line(0,1){16}} \put(301,55){\line(1,0){23}}
\put(300,12){\line(0,1){43}} \put(373,12){\line(-1,0){73}}
\put(279,169){\line(1,0){92}} \put(279,130){\line(0,1){39}}
\put(235,130){\line(1,0){43}} \put(234,130){\line(0,-1){73}}
\put(75,31){\line(1,0){102}} \put(177,165){\line(0,-1){133}}
\put(39,167){\line(1,0){138}} \put(75,63){\line(0,-1){31}}
\put(62,64){\line(1,0){13}} \put(60,79){\line(0,-1){15}}
\put(99,79){\line(-1,0){39}} \put(99,102){\line(0,-1){23}}
\put(134,102){\line(-1,0){34}} \put(134,155){\line(0,-1){52}}
\put(68,156){\line(1,0){66}} \put(68,139){\line(0,1){17}}
\put(38,139){\line(0,1){28}} \put(38,139){\line(1,0){30}} \thinlines   
\put(12,62){\line(1,0){85}} \put(12,62){\line(0,1){74}}
\put(96,136){\line(0,-1){74}} \put(12,136){\line(1,0){84}} \end{picture}}

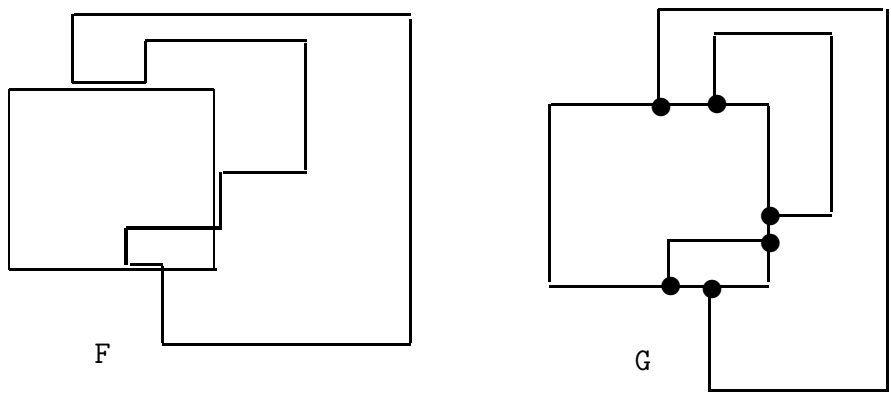
\begin{center} {\bf Figure 3 -- Second Example of Coloring and Formation}
\end{center} \vspace{3mm}

\section{\bf Simple Operations and the Parity Lemma}

Recall that a {\em circuit} in a graph $G$ is a subgraph that is equivalent to a
circle graph (i.e. homeomorphic to a circle). \vspace{3mm}

Let $G$ be a cubic graph. Suppose that $C$ is a coloring of $G$ with three colors
(so that three distinct colors are incident at each vertex of $G$). Let the
colors be denoted by $r$ (red), $b$ (blue) and $p$ (purple). Then we can classify
circuits in $G$ relative to the coloring $C$.  We shall be concerned with those
circuits that contain exactly two colors. The possible two-color circuits are
$r-b$ (red-blue), $r-p$ (red-purple) and $b-p$ (blue-purple). Let $\Delta(G,C)$
denote the number of distinct two-color circuits in $G$ with the coloring $C$.
\vspace{2mm}

\noindent {\bf Definition.} Call the {\em parity} of the coloring $C$, denoted
$\pi(G,C)$, the parity of the number of distinct two-color circuits,
$\Delta(G,C)$. \vspace{3mm}

\noindent {\bf Definition.}  If $C$ is a coloring of $G$ and $d$ is a two-color
circuit in $G$ (called a {\em modulus} in \cite{SB}), then we can obtain a new coloring $C' = C$ of $G$
by interchanging the colors on $d$.  Call the operation of switching colors on a two
color circuit a {\em simple operation} on the coloring $C$. \vspace{3mm}

In this section we will prove a basic Parity Lemma due to  Spencer-Brown
\cite{SB} in the category of formations.  A similar result due to W. T. Tutte \cite{T} in
the category of plane cubic graphs implies the Parity Lemma, but is proved by
a different method.  The lemma states that simple operations
on planar graphs or planar formations preserve parity.  Note that  by the results
of section one, colorings of cubic graphs and formations for cubic graphs are in
one-to-one correspondence. The proof of the parity lemma given here is due to the
author of this paper. \vspace{3mm}

Note that for a formation $F$ composed of red and blue curves, the two color
circuits are counted by $\Delta(F) = R + B + Alt$ where $R$ denotes the number of
red curves, $B$ denotes the number of blue curves, and $Alt$ denotes the number of
red-blue alternating circuits in the corresponding coloring. These red-blue
circuits are characterized in the formation as those two-colored circuits that avoid the
places where there is a superposition of red and blue (these places correspond to
purple edges in the coloring). The red curves in the formation correspond to
red-purple circuits in the coloring, and the blue curves in the formation
correspond to blue-purple circuits in the coloring. \vspace{3mm} 

Each formation corresponds to a specific graph coloring. Simple operations on the coloring
induce new formations over the underlying graph. Simple operations can be
performed directly on a formation via a graphical calculus. This calculus is
based on the principle of {\em idemposition} saying that: {\em superposition of segments
of the same color results in the cancellation of those segments}. The result of
an idemposition of curves of the same color is a mod-$2$ addition of the curves.
Two curves of the same color that share a segment are joined at the junctions of
the segment, and the segment disappears. In order to perform a simple operation
on a blue loop, superimpose a red loop upon it and perform the corresponding
idemposition with the other red curves that impinge on this red loop along the
blue loop. Similarly, in order to perform a simple operation on a red loop,
superimpose a blue loop on it and idempose this blue loop with the blue curves
that impinge on the red loop. Finally, in order to perform a simple operation on
a red-blue alternating circuit in a formation, superimpose a red and a blue loop on this
circuit and perform the corresponding idempositions. These instructions for
performing simple operations are illustrated in Figure 4. In this Figure some of the edges that 
are intended to be superimposed are drawn at a short distance from one another in order to enhance
the reader's ability to trace the curves.
\vspace{3mm}

{\tt    \setlength{\unitlength}{0.92pt} \begin{picture}(358,503) \thinlines   
\put(139,29){\makebox(192,26){idempose blue and red curves}}
\put(125,132){\makebox(98,22){operate on alt}} \put(214,89){\vector(-3,-1){67}}
\put(143,126){\vector(1,0){51}} \put(68,10){\framebox(53,64){}} \thicklines  
\put(21,14){\framebox(96,56){}} \thinlines    \put(240,106){\framebox(26,51){}}
\thicklines   \put(235,102){\framebox(37,59){}} \thinlines   
\put(232,97){\framebox(97,69){}} \thicklines   \put(276,102){\framebox(49,60){}}
\thinlines    \put(19,99){\framebox(97,69){}} \thicklines  
\put(63,104){\framebox(49,60){}} \thinlines   
\put(174,209){\makebox(131,29){idempose blue curves}}
\put(120,291){\makebox(95,27){operate on red}} \put(207,253){\vector(-3,-1){67}}
\put(144,278){\vector(1,0){51}} \put(107,326){\line(0,-1){39}}
\put(101,216){\line(0,-1){15}} \put(112,249){\line(0,-1){32}}
\put(13,186){\line(1,0){117}} \put(13,216){\line(0,-1){31}}
\put(129,217){\line(0,-1){32}} \put(113,217){\line(1,0){16}}
\put(31,216){\line(-1,0){18}} \put(108,195){\line(0,1){1}}
\put(86,216){\line(1,0){15}} \put(49,216){\line(1,0){21}} \thicklines  
\put(15,187){\framebox(112,27){}} \thinlines    \put(49,216){\line(0,1){20}}
\put(49,236){\line(1,0){36}} \put(86,236){\line(0,-1){20}}
\put(70,216){\line(0,-1){15}} \put(70,201){\line(1,0){31}}
\put(112,249){\line(-1,0){81}} \put(31,249){\line(0,-1){33}}
\put(240,266){\framebox(103,19){}} \thicklines  
\put(236,262){\framebox(112,27){}} \thinlines    \put(251,292){\line(1,0){18}}
\put(269,292){\line(0,1){20}} \put(269,312){\line(1,0){36}}
\put(306,312){\line(0,-1){20}} \put(306,292){\line(-1,0){16}}
\put(290,292){\line(0,-1){15}} \put(290,277){\line(1,0){31}}
\put(321,278){\line(0,1){7}} \put(321,286){\line(1,0){12}}
\put(333,288){\line(0,1){37}} \put(332,325){\line(-1,0){81}}
\put(251,325){\line(0,-1){33}} \put(25,326){\line(0,-1){33}}
\put(106,326){\line(-1,0){81}} \put(95,287){\line(1,0){12}}
\put(95,279){\line(0,1){7}} \put(64,278){\line(1,0){31}}
\put(64,293){\line(0,-1){15}} \put(80,293){\line(-1,0){16}}
\put(80,313){\line(0,-1){20}} \put(43,313){\line(1,0){36}}
\put(43,293){\line(0,1){20}} \put(25,293){\line(1,0){18}} \thicklines  
\put(10,264){\framebox(112,27){}} \thinlines   
\put(158,367){\makebox(113,36){idempose red curves}}
\put(201,418){\vector(-3,-1){67}} \put(138,442){\vector(1,0){51}} \thicklines  
\put(129,439){\makebox(82,30){operate on blue}}
\put(223,417){\framebox(104,76){}} \thinlines   
\put(226,421){\framebox(97,69){}} \thicklines   \put(270,426){\framebox(49,60){}}
\thinlines    \put(18,344){\framebox(104,66){}} \thicklines  
\put(22,347){\framebox(47,59){}} \put(62,427){\framebox(49,60){}} \thinlines   
\put(18,422){\framebox(97,69){}} \end{picture}}

\begin{center} {\bf Figure 4  -- Simple Operations} \end{center} \vspace{20mm}

\noindent {\bf Parity Lemma}. If $C'$ and $C$ are colorings of a planar cubic
graph $G$ with  $C'$ obtained from $C$ by a simple operation, then the parity of
$C'$ is equal to the parity of $C$,  $\pi(C') = \pi(C).$  Equivalently, parity is
preserved under simple operations on planar formations. \vspace{3mm}

In order to prove the Parity Lemma, we need to consider elementary properties of
idempositions of curves in the plane. \vspace{3mm}

First consider the idemposition of two curves of the same color, as illustrated in
Figure 5. We can distinguish three types of interaction denoted by $L$
(Left), $R$ (Right) and $B$ (bounce).  A bounce ($B$) is when the second curve
shares a segment with the first curve, but does not cross the first curve. 
Crossing interactions are classified as left and right according as the person
walking along the first curve first encounters the second curve on his right
($R$) or on his left ($L$). After an encounter, there ensues a shared segment
that the walker leaves in the direction of the opposite hand.  Let $|L|$ denote the
number of left crossings between the first and second curves, $|R|$ the number of
right crossings, and $|B|$ the number of bounces. (Note the $|L|$ and $|R|$
depend upon the choice of direction for the walk along the first curve.  Let
$P(A,A')$ denote the parity of $(|L|-|R|)/2 + |B|$ for an interaction of curves
$A$ and $A'.$ \vspace{20mm}

{\tt    \setlength{\unitlength}{0.92pt} 
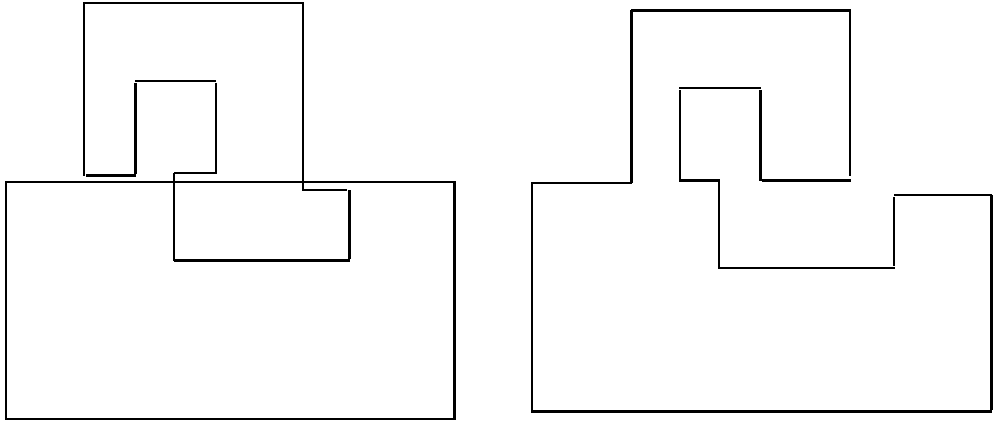
\begin{picture}(427,242) \thinlines   
\put(106,10){\makebox(181,35){(|L|-|R|)/2 + |B| = (1-1)/2 +1 = 1}}
\put(226,62){\line(1,0){189}} \put(415,151){\line(0,-1){88}}
\put(357,227){\line(0,-1){68}} \put(226,156){\line(0,-1){94}}
\put(375,151){\line(1,0){40}} \put(321,157){\line(1,0){36}}
\put(287,157){\line(1,0){16}} \put(226,156){\line(1,0){41}}
\put(287,157){\line(0,1){37}} \put(287,195){\line(1,0){33}}
\put(320,194){\line(0,-1){37}} \put(303,157){\line(0,-1){36}}
\put(303,121){\line(1,0){72}} \put(375,122){\line(0,1){28}}
\put(267,227){\line(0,-1){71}} \put(267,227){\line(1,0){90}}
\put(42,230){\line(1,0){90}} \put(42,230){\line(0,-1){71}}
\put(132,153){\line(0,1){77}} \put(150,153){\line(-1,0){18}}
\put(151,125){\line(0,1){28}} \put(79,124){\line(1,0){72}}
\put(79,160){\line(0,-1){36}} \put(96,160){\line(-1,0){17}}
\put(96,197){\line(0,-1){37}} \put(63,198){\line(1,0){33}}
\put(63,160){\line(0,1){37}} \put(43,159){\line(1,0){20}}
\put(10,59){\framebox(184,97){}} \end{picture}}

\begin{center} {\bf Figure 5 -- Idemposition} \end{center} \vspace{3mm}

\noindent {\bf Idemposition Lemma.}  Let  $A$ and $A'$ be two simple closed
curves in the plane of the same color.  The parity of the number of simple closed
curves resulting from the idemposition of $A$ and $A'$ is equal to $P(A,A') = 
(|L|-|R|)/2 + |B|$ (mod $2$) where the terms in this formula are  as defined
above. \vspace{3mm}

\noindent {\bf Proof.}  The proof is by induction on the number of crossing
interactions between the two curves.  It is easy to see that the removal of a
bounce changes the parity of the idemposition. See Figure 6.  Figure 7
illustrates a collection of unavoidable crossing  interactions between two
curves.  That is, if there are crossing interactions, then one of the situations
in Figure 7 must occur. ( To see this note that if you follow curve $A'$ and
cross curve $A$, then there is a first place where $A'$ crosses $A$ again.  The
unavoidable configurations are a list of the patterns of crossing and crossing
again.) It is then clear from Figure 7, by counting parity after the indicated
idempositions, that the result follows by induction. // \vspace{3mm}

{\tt    \setlength{\unitlength}{0.92pt} \begin{picture}(173,87) \thinlines   
\put(10,12){\line(1,0){153}} \put(140,44){\line(0,0){0}}
\put(101,17){\line(0,1){59}} \put(65,16){\line(1,0){36}}
\put(65,77){\line(0,-1){61}} \end{picture}}

\begin{center} {\bf Figure 6 -- Bounce} \end{center} \vspace{10mm}

{\tt    \setlength{\unitlength}{0.92pt} \begin{picture}(285,178) \thinlines   
\put(231,35){\line(0,-1){22}} \put(232,35){\line(1,0){10}}
\put(242,35){\line(0,1){41}} \put(186,12){\line(1,0){45}}
\put(186,40){\line(0,-1){28}} \put(203,41){\line(-1,0){17}}
\put(204,74){\line(0,-1){32}} \put(76,43){\line(0,1){34}}
\put(90,42){\line(-1,0){14}} \put(91,14){\line(0,1){28}}
\put(39,13){\line(1,0){52}} \put(39,43){\line(0,-1){30}}
\put(56,43){\line(-1,0){17}} \put(56,79){\line(0,-1){36}}
\put(216,123){\line(0,1){43}} \put(231,123){\line(-1,0){14}}
\put(231,91){\line(0,1){31}} \put(200,90){\line(1,0){31}}
\put(199,122){\line(0,-1){32}} \put(183,122){\line(1,0){16}}
\put(183,167){\line(0,-1){45}} \put(99,123){\line(0,1){45}}
\put(79,123){\line(1,0){20}} \put(79,94){\line(0,1){29}}
\put(49,93){\line(1,0){30}} \put(49,123){\line(0,-1){30}}
\put(29,123){\line(1,0){20}} \put(29,168){\line(0,-1){45}}
\put(159,38){\line(1,0){116}} \put(10,119){\line(1,0){116}}
\put(154,120){\line(1,0){116}} \put(12,39){\line(1,0){116}} \end{picture}}

\begin{center} {\bf Figure 7 --  Innermost Cross and Recross} \end{center}
\vspace{3mm}

\noindent {\bf Proof of the Parity Lemma.} Consider a formation $F$ consisting of
one red loop $A$ that is touched by a set of  $n$ disjoint blue curves. It is 
clear by construction that the number of alternating (red/blue) circuits in $F$
is equal to the number of curves in the idemposition obtained after letting all
the blue curves become red (so that they cancel with the original red loop where
blue meets red).  As a result, we can apply the Idemposition Lemma to conclude
that $$\Delta(F) =  1 + n + (|L|-|R|)/2 + |B| (mod 2)$$ where $|L|$,$|R|$ and
$|B|$ denote the total number of left, right and bounce interactions between the
red curve and the blue curves and $n$ is the number of blue curves. (Apply the
Lemma to each blue curve one at a time.) The main point is that the parity of $F$
is determined by a count of local interactions along the red curve $A$. In this
case, when we perform a simple operation on $A$, the curve count does not change.
We simply interchange the roles of blue circuits (i.e. blue/purple circuits) and
alternating circuits (i.e. red/blue circuits).  Thus in this case we have that
$\Delta(F) = \Delta(F')$ where $F'$ is obtained by a simple operation on the
curve $A$ in $F$.  Hence parity is certainly preserved. \vspace{3mm}

In the general case we have a red curve $A$ that interacts with a collection of
blue curves, and these blue curves interact with the rest of the formation. Call
the whole formation $F$, and let $G$ denote the subformation consisting of the
curve $A$ and all the blue curves that interact with $A$.  If $F'$ is the result
of operating on $A$ in $F$, then $F'$ will contain $G'$, the result of operating
on $A$ in $G$.  $G'$ will consist of the curve $A$ plus all blue curves in $F'$
that touch the curve $A$ in $F'$.  In counting the change of $\Delta$ from
$\Delta(F)$ to $\Delta(F')$, we actually count the change in the count of blue
curves and the change in the count of alternating circuits.  Each of these
changes can be regarded as the result of a single color idemposition originating
at $A$.  The change in $\Delta$ from $F$ to $F'$ is the sum of the change in the number of blue
curves and the change in the number of  alternating circuits. Each of these changes is
determined by local interactions along the curve $A$. The parity of the change again
depends only on these local interactions. Since the transformation from $G$ to
$G'$ has identical local interactions, and since $G$ and $G'$ have the same
$\Delta$ and hence the same parity, it follows that $F$ and $F'$ have the same
parity.  This completes the proof of the Parity Lemma. //
\vspace{3mm}

\noindent {\bf Remark.} In performing a simple operation, the curve count may change without 
changing the parity of this count. Figure 8 illustrates an example of this phenomenon.
\vspace{20mm}

{\tt    \setlength{\unitlength}{0.92pt} \begin{picture}(317,168) \thicklines  
\put(10,12){\framebox(135,85){}} \thinlines    \put(21,101){\line(1,0){21}}
\put(42,99){\line(0,-1){22}} \put(42,77){\line(1,0){30}}
\put(72,78){\line(0,1){21}} \put(72,100){\line(-1,0){12}}
\put(60,100){\line(0,1){50}} \put(21,151){\line(1,0){39}}
\put(95,153){\line(0,-1){54}} \put(95,99){\line(1,0){27}}
\put(122,99){\line(0,1){54}} \put(95,153){\line(1,0){27}}
\put(21,151){\line(0,-1){50}} \thicklines   \put(62,118){\framebox(31,22){}}
\put(227,121){\framebox(31,22){}} \thinlines    \put(185,154){\line(0,-1){50}}
\put(259,156){\line(1,0){27}} \put(259,156){\line(0,-1){54}}
\put(185,154){\line(1,0){39}} \put(224,103){\line(0,1){50}}
\put(236,81){\line(0,1){21}} \put(206,80){\line(1,0){30}}
\put(206,102){\line(0,-1){22}} \thicklines   \put(167,15){\framebox(135,85){}}
\thinlines    \put(206,102){\line(1,0){18}} \put(236,102){\line(1,0){23}}
\put(185,103){\line(-1,0){22}} \put(162,103){\line(0,-1){91}}
\put(162,12){\line(1,0){143}} \put(305,12){\line(0,1){92}}
\put(305,104){\line(-1,0){19}} \put(286,156){\line(0,-1){52}} \end{picture}}

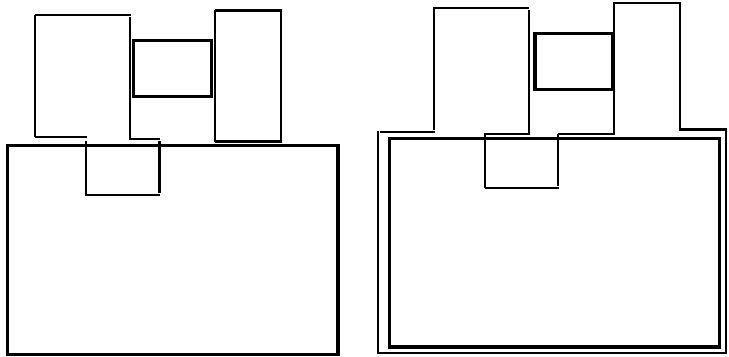
\begin{center} {\bf Figure 8 --  Changing Curve Count Under Simple Operation}
\end{center} \vspace{3mm}

{\tt    \setlength{\unitlength}{0.92pt} \begin{picture}(317,168) \thinlines   
\put(286,156){\line(0,-1){52}} \put(305,104){\line(-1,0){19}}
\put(305,12){\line(0,1){92}} \put(162,12){\line(1,0){143}}
\put(162,103){\line(0,-1){91}} \put(185,103){\line(-1,0){22}}
\put(236,102){\line(1,0){23}} \put(206,102){\line(1,0){18}} \thicklines  
\put(167,15){\framebox(135,85){}} \thinlines    \put(206,102){\line(0,-1){22}}
\put(206,80){\line(1,0){30}} \put(236,81){\line(0,1){21}}
\put(224,103){\line(0,1){50}} \put(185,154){\line(1,0){39}}
\put(259,156){\line(0,-1){54}} \put(259,156){\line(1,0){27}}
\put(185,154){\line(0,-1){50}} \put(21,151){\line(0,-1){50}}
\put(95,153){\line(1,0){27}} \put(122,99){\line(0,1){54}}
\put(95,99){\line(1,0){27}} \put(95,153){\line(0,-1){54}}
\put(21,151){\line(1,0){39}} \put(60,100){\line(0,1){50}}
\put(72,100){\line(-1,0){12}} \put(72,78){\line(0,1){21}}
\put(42,77){\line(1,0){30}} \put(42,99){\line(0,-1){22}}
\put(21,101){\line(1,0){21}} \thicklines   \put(10,12){\framebox(135,85){}}
\end{picture}}

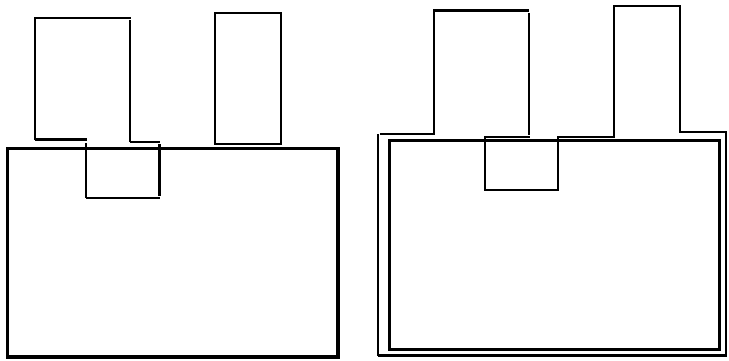
\begin{center} {\bf Figure 9 --  Unchanging Curve Count Under Simple Operation}
\end{center} \vspace{3mm}

\noindent {\bf Remark.}  The parity lemma fails for a non-planar formation.
For example, consider the formation in Figure 10.  This is a formation for the
Petersen graph with one edge removed.  As the Figure indicates, parity is not
preserved by a simple operation on this graph. The curve count in the first
formation is five and the curve count in the second formation (after the simple
operation) is four. This shows that the underlying graph of these two formations 
is non-planar. \vspace{3mm}

{\tt    \setlength{\unitlength}{0.92pt} \begin{picture}(360,163) \thinlines   
\put(197,106){\framebox(147,40){}} \thicklines   \put(348,150){\line(0,-1){47}}
\put(192,150){\line(0,-1){48}} \put(348,151){\line(-1,0){156}}
\put(323,103){\line(1,0){25}} \put(266,103){\line(1,0){14}}
\put(229,102){\line(-1,0){37}} \put(281,50){\line(1,0){23}}
\put(240,50){\line(1,0){27}} \put(265,102){\line(3,-4){39}}
\put(282,50){\line(4,5){42}} \put(229,102){\line(3,-4){39}}
\put(239,49){\line(4,5){42}} \thinlines    \put(196,10){\framebox(152,37){}}
\put(12,14){\framebox(152,37){}} \put(10,111){\framebox(152,37){}} \thicklines  
\put(47,54){\line(4,5){42}} \put(37,107){\line(3,-4){39}}
\put(90,55){\line(4,5){42}} \put(73,107){\line(3,-4){39}}
\put(37,108){\line(1,0){36}} \put(88,108){\line(1,0){43}}
\put(48,54){\line(1,0){27}} \put(90,55){\line(1,0){23}} \end{picture}}

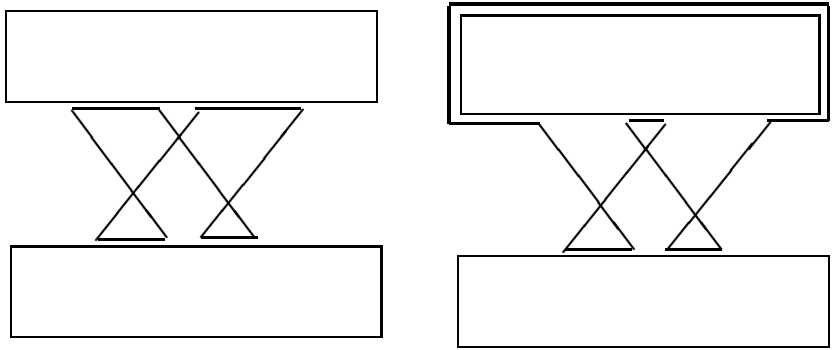
\begin{center} {\bf Figure 10 --  Parity Reversed in the One-deleted Petersen}
\end{center} \vspace{3mm}

\section{\bf A Principle of Irreducibility}

The main result of this section is the equivalence of the Four Color Theorem with
a property of formations that I call the {\em Primality Principle}. In order to
state this property we need to explain the concept of a {\em trail} in a
formation, and {\em how a trail can facilitate or block an  attempt to extend a
coloring.} \vspace{3mm}

Consider a formation with two blue curves and a single red curve that interacts
with the two blues.  See Figure 11 for an illustration of this condition. I shall
call the red curve a {\em trail} between the two blues.  Call the blue curves the
{\em containers} or {\em contextual curves} for the trail. Call the {\em graph of the trail $T$}
the cubic
graph $G(T)$  corresponding to the formation consisting in the two blues and the
red curve between them as in Figure 12. In Figure 11 we have also indicated a
double arrow pointing between the two blue curves and disjoint from the trail. The
double arrow is meant to indicate an edge that we would like to color, extending the
given formation to a new formation that includes this edge. We shall refer to this 
double arrow as the {\em empty edge}.  In the example shown
in Figure 13 we obtain this extension by drawing a purple (blue plus red)
curve that goes through the empty edge. The part of the purple curve that is not
on the arrow forms a pathway in the given formation from  one arrow-tip to the other
that  uses only two colors (red and blue). After idemposition, this purple curve
effects a two-color switch along this pathway and the formation is extended as
desired. Under these circumstances we say that the formation is {\em completable over the 
empty edge.} If simple operations on a given formation with an empty edge can transform it so that
the formation is completable over the empty edge, we say that the formation is {\em completable by simple operations.}
Since the final action of completing the formation changes the empty edge to a colored edge, this last operation
(described above) will be called a {\em complex operation.} \vspace{3mm}

{\tt    \setlength{\unitlength}{0.92pt} \begin{picture}(223,217) \thicklines  
\put(46,110){\vector(0,-1){22}} \put(46,91){\vector(0,1){39}}
\put(120,84){\line(1,0){35}} \put(155,134){\line(0,-1){49}}
\put(120,135){\line(0,-1){51}} \put(106,135){\line(1,0){14}}
\put(168,135){\line(-1,0){13}} \put(169,159){\line(0,-1){23}}
\put(107,159){\line(1,0){62}} \put(106,135){\line(0,1){24}} \thinlines   
\put(10,10){\framebox(203,76){}} \put(10,131){\framebox(203,76){}}
\put(161,98){\makebox(31,23){T}} \end{picture}}

\begin{center} {\bf Figure 11 --  A Trail Between Two Blue curves} \end{center}
\vspace{3mm}

{\tt    \setlength{\unitlength}{0.92pt} \begin{picture}(438,218) \thinlines   
\put(308,33){\makebox(36,27){$G^{*}(T)$}} \put(94,30){\makebox(39,31){$G(T)$}}
\put(370,87){\circle*{8}} \put(335,85){\circle*{8}} \put(369,133){\circle*{8}}
\put(335,133){\circle*{8}} \put(382,133){\circle*{8}} \put(321,133){\circle*{8}}
\put(258,86){\circle*{8}} \put(258,133){\circle*{8}} \put(155,87){\circle*{8}}
\put(121,86){\circle*{8}} \put(168,133){\circle*{8}} \put(155,133){\circle*{8}}
\put(121,133){\circle*{8}} \put(107,134){\circle*{8}}
\put(225,132){\framebox(203,76){}} \put(224,10){\framebox(203,76){}}
\put(258,132){\line(0,-1){46}} \put(335,132){\line(0,-1){46}}
\put(369,132){\line(0,-1){46}} \put(321,159){\line(1,0){61}}
\put(321,159){\line(0,-1){26}} \put(382,159){\line(0,-1){26}}
\put(168,159){\line(0,-1){26}} \put(107,159){\line(0,-1){26}}
\put(107,159){\line(1,0){61}} \put(155,132){\line(0,-1){46}}
\put(121,132){\line(0,-1){46}} \put(10,10){\framebox(203,76){}}
\put(11,132){\framebox(203,76){}} \end{picture}}

\begin{center} {\bf Figure 12 --  The Graphs $G(T)$ and $G^{*}(T)$} \end{center}
\vspace{3mm}

{\tt    \setlength{\unitlength}{0.92pt} \begin{picture}(280,173) \thinlines   
\put(189,73){\framebox(22,26){}} \thicklines   \put(185,70){\framebox(30,33){}}
\thinlines    \put(163,104){\framebox(106,58){}}
\put(164,10){\framebox(106,58){}} \thicklines   \put(217,64){\framebox(36,44){}}
\put(34,85){\vector(0,1){18}} \put(34,103){\vector(0,-1){31}}
\put(64,65){\framebox(36,44){}} \thinlines    \put(11,11){\framebox(106,58){}}
\put(10,105){\framebox(106,58){}} \end{picture}}

\begin{center} {\bf Figure 13 -- A Colorable Trail} \end{center} \vspace{3mm}

{\tt    \setlength{\unitlength}{0.92pt} \begin{picture}(173,152) \thinlines   
\put(10,10){\framebox(152,37){}} \put(11,105){\framebox(152,37){}} \thicklines  
\put(67,50){\line(4,5){42}} \put(57,102){\line(3,-4){39}}
\put(108,50){\line(4,5){42}} \put(93,102){\line(3,-4){39}}
\put(57,103){\line(1,0){36}} \put(107,103){\line(1,0){43}}
\put(68,49){\line(1,0){27}} \put(110,50){\line(1,0){23}}
\put(32,103){\vector(0,-1){52}} \put(32,73){\vector(0,1){30}} \end{picture}}

\begin{center} {\bf Figure 14 -- The Petersen Trail} \end{center} \vspace{3mm}

Another example of a trail is shown in Figure 14.  Here no extension is possible
since the extended graph is the Petersen graph, a graph that does not admit a
coloration.
\bigbreak

In a trail the endpoints of the empty edge are specified, since one would
like to complete the formation over the empty edge. The simplest example of
uncolorability is just two curves and an empty edge. Then no matter how the
curves are colored there is no way to extend the formation over the empty
edge.
\bigbreak

There are two cases in the coloring structure of a trail: The two contextual curves have the same color or they
have different colors. We shall distinguish these two cases by {\em defining} those colors of 
the contextual curves to be
the colors incident at the endpoints of the empty edge. Note that when we refer to a curve in
a formation we mean either a blue curve, a red curve or a cycle that alternates in red and blue when we
are performing a parity count. On the other hand one can also consider purple curves, but these will
appear in the formation as alternations of purple with blue or red at those sites where
the purple is idemposed
with red or blue respectively. When counting curves, we shall only count blue, red
and alternating (red and blue).
\bigbreak

We first consider contextual curves of the same color.
Suppose that both contextual curves are purple. Then any
trail between them must be drawn either in blue or in red. It will be called
a {\em non-purple trail.} Thus one could insert a Petersen trail drawn as a
red curve and then idemposed between the purple curves, or a Petersen trail drawn
as a blue curve and then idemposed between the two purples. We will say that a trail between 
two curves of the same color (red,blue or purple) is {\em factored} if after removing
the two contextual curves (by idemposing them with curves of the same color) the remaining trail
structure has multiple components. 
\bigbreak

See Figure 15 for an illustration of this removal process. 
The ``trail" that we 
uncover by the removal process is {\em not} the color of the two curves and it does not touch
the endpoints of the empty edge. In Figure 15 we illustrate a trail between two purples. That is, 
each endpoint of the empty edge touches the color purple.  In the second part of the figure we reveal the purples so 
that an idemposition of this figure gives the first part and a removal of the two purple curves gives
the single trail component. Since there is only one component, this trail is not factored. 
\bigbreak

{\tt    \setlength{\unitlength}{0.92pt}
\begin{picture}(292,295)
\thicklines   \put(46,25){\vector(0,1){57}}
              \put(46,83){\vector(0,-1){63}}
              \put(42,190){\vector(0,1){46}}
              \put(42,237){\vector(0,-1){59}}
              \put(182,247){\line(1,0){103}}
              \put(78,245){\line(-1,0){73}}
              \put(104,244){\line(1,0){50}}
              \put(156,163){\line(1,0){134}}
              \put(103,164){\line(-1,0){102}}
\thinlines    \put(5,252){\line(1,0){278}}
              \put(3,171){\line(1,0){285}}
\thicklines   \put(78,246){\line(0,1){46}}
              \put(79,292){\line(1,0){102}}
              \put(181,292){\line(0,-1){45}}
              \put(156,244){\line(0,-1){81}}
              \put(104,164){\line(0,1){80}}
              \put(108,3){\line(0,1){80}}
              \put(157,3){\line(-1,0){48}}
              \put(157,84){\line(0,-1){81}}
              \put(182,85){\line(-1,0){25}}
              \put(182,130){\line(0,-1){45}}
              \put(80,130){\line(1,0){102}}
              \put(80,84){\line(0,1){46}}
              \put(108,84){\line(-1,0){26}}
\thinlines    \put(6,13){\line(1,0){285}}
              \put(8,94){\line(1,0){278}}
\thicklines   \put(6,10){\line(1,0){282}}
              \put(6,90){\line(1,0){280}}
\end{picture}}

\begin{center} {\bf Figure 15 --- A Trail Between Two Purples} \end{center}
\bigbreak

 For the case where both contextual curves have the same color there is no loss of generality
in assuming that the two contextual curves are both red or both blue. Then any extra curves produced
in a factorization can be seen directly, in their appearance as alternating, blue or red.
\bigbreak

Secondly, suppose that the two contextual curves have different colors. And suppose that the
formation has a non-empty trail structure between the two curves. We say that this
formation is {\em factored} with respect to the empty edge if there is an extra curve
in the formation that does not pass through either endpoint of the empty
edge. For example, perform a simple operation on the Petersen trail as in Figure 14, making the top
curve purple (at an endpoint of the empty edge). Note that in this example
every curve in the formation
passes through one of the endpoints of the empty edge, so it is not factored. 
Second example: operate on the top curve in
Figure 11. You will find that this produces a red curve that is isolated
from the empty edge, giving a factorization.
\bigbreak 

We shall say that a formation is {\em unfactored} if it is not factored.
\bigbreak 

We shall say that a trail $T$ {\em factorizes} if there is a
formation for the graph $G(T)$ (see definition above) of this trail that is factored. 
Note that we do not require that the original version of the trail be factored. A new version 
can be obtained by simple operations on the original formation, or by more complicated re-colorings.
A trail is said to be {\em prime} if it does not admit
any factorization.
\bigbreak

Sometimes a trail can factorize by simple operations as in the example in Figure
16.  In the example in Figure 16 we perform a simple operation on the upper blue
curve. Note that in the resulting factorization the arrow is now between a lower
blue curve and a part of the upper blue curve that has a superimposed red segment
from one of the factors. The Petersen trail of Figure 14 is a significant example of a
prime trail. Recolorings of the graph of the formation of this trail just return the
Petersen trail in slightly disguised form. \bigbreak

{\tt    \setlength{\unitlength}{0.92pt} \begin{picture}(453,234) \thicklines  
\put(441,222){\line(0,-1){81}} \put(405,141){\line(1,0){34}}
\put(229,222){\line(1,0){211}} \put(228,141){\line(0,1){81}}
\put(340,141){\line(-1,0){111}} \put(355,141){\line(1,0){36}}
\put(356,89){\line(1,0){35}} \thinlines    \put(232,10){\framebox(203,76){}}
\put(234,143){\framebox(203,76){}} \thicklines   \put(340,141){\line(0,1){24}}
\put(341,166){\line(1,0){62}} \put(405,166){\line(0,-1){23}}
\put(355,141){\line(0,-1){51}} \put(391,139){\line(0,-1){49}}
\put(263,98){\vector(0,1){39}} \put(263,116){\vector(0,-1){22}}
\put(47,114){\vector(0,-1){22}} \put(47,95){\vector(0,1){39}}
\put(121,88){\line(1,0){35}} \put(156,138){\line(0,-1){49}}
\put(121,139){\line(0,-1){51}} \put(107,139){\line(1,0){14}}
\put(169,139){\line(-1,0){13}} \put(170,163){\line(0,-1){23}}
\put(108,163){\line(1,0){62}} \put(107,139){\line(0,1){24}} \thinlines   
\put(11,14){\framebox(203,76){}} \put(10,136){\framebox(203,76){}}
\put(162,102){\makebox(31,23){T}} \end{picture}}

\begin{center} {\bf Figure 16 -- A Factorizable Trail} \end{center} \vspace{3mm}

A trail is said to be {\em uncolorable} if the graph $G^{*}(T)$ obtained from
$G(T)$ by adding the edge corresponding to the double arrow is an uncolorable
graph.  Thus the Petersen trail is uncolorable since $G^{*}(T)$ is the Petersen
graph. A trail is said to be a {\em minimal uncolorable} trail if the graph
$G^{*}(T)$ is a smallest uncolorable graph. Now the Petersen graph is the
smallest possible uncolorable graph other than the dumbell shown in Figure 1.  In
particular, the Petersen is the smallest non-planar uncolorable.  This does not,
in itself, rule out the possibility of planar uncolorables other than the dumbell
(that is the essence of the Four Color Theorem).  Hence we can entertain the {\em
possibility} of minimal planar uncolorable trails. \vspace{3mm}

Now we can state the \vspace{3mm}

\noindent {\bf Primality Principle.}  A minimal planar (non-empty) uncolorable
trail is prime. \vspace{3mm}

In other words, this principle states that there is no possibility of making a
minimal planar uncolorable trail that is factored into smaller planar trails. The
principle lends itself to independent investigation since one can try combining
trails to make a possibly uncolorable formation (i.e. that the graph
$G^{*}(T_{1}, ... , T_{n})$ is uncolorable where this graph is obtained from the
formation consisting in the trails $T_{1}, ... , T_{n}$ placed disjointly between
two blue curves.) The combinatorics behind this principle are the subject of much
of the research of G. Spencer-Brown.  Spencer-Brown regards the Primality
Principle as {\em axiomatic} (See \cite{SB} page 169).  It  is one purpose of
this paper to point out the equivalence of the Four Color Theorem and the Primality
Principle.  
\vspace{3mm}

\noindent {\bf Theorem.}  The Primality Principle is {\em equivalent} to the Four
Color Theorem. \vspace{2mm}

\noindent {\bf Proof.} First suppose the Primality Principle - that minimal
uncolorable trails are prime. Let $T$ be a minimal uncolorable non-empty planar
trail.  Without loss of generality $T$ is defined by a formation consisting in a single red
curve (the trail) drawn between two disjoint blue curves. We call this formation $F(T)$, the
formation induced by the trail $T$. The formation can be depicted so that the two blue curves
appear as parallel lines (to be completed to circuits -- above for the top line and
below for the bottom line) and the trail $T$ is interacting between the two parallel
blue lines. In this depiction, we can set a double arrow indicator between the
two parallel lines, with this indicator entirely to the left of $T$. This double-arrow indicator
represents
an edge that we would like to complete to form a larger formation/coloring. 
Uncolorability of the trail means that there is no coloring of the graph obtained
by adding to the underlying graph of $F(T)$ an edge corresponding to the double arrow. \bigbreak

Note that an uncolorable trail is necessarily incompletable (across the empty edge) by
simple operations. This implies that there is no two-color pathway in the given formation
of the trail from one endpoint of the empty edge to the other endpoint. We can use these facts to count
the number of curves in a minimal uncolorable trail.
\bigbreak

First consider a prime uncolorable trail
with two blue  contextual curves, and an existing trail between them in red. This trail must consist in
a single red curve. There can be no other red curves in the formation. There
is an alternating curve incident to each endpoint of the empty edge. Thus there are at most two
alternating curves, one for each endpoint of the empty edge. (Other alternating curves would become 
red components after the removal of the contextual curves.) If there is one alternating curve, then
there is a two-color pathway between the endpoints of the empty edge, and the formation is 
completable over this edge. Hence there are two alternating curves.
Thus we see that the curve count (one red, two blue, two alternating) for a prime
uncolorable formation with two blue contextual curves is {\em five.}
\bigbreak

Second, consider a prime, uncolorable trail with one purple contextual curve
and one blue contextual curve. The trail structure will then consist in red curves woven 
between the two 
contextual curves. Once these red curves are idemposed with the purple, the formation can be
regarded as two blue curves with the trail structure passing through (say) the 
upper endpoint of the empty edge, so that this upper endpoint rests on purple. 
Such a formation is unfactored if and only if all curves  pass through the endpoints of the empty 
edge. Thus we have a single blue curve and a single alternating curve passing through the lower 
endpoint, and one red curve and one blue curve passing through the upper endpoint.  This makes a total
of two blues, one red and one alternator, hence a curve count of {\em four}\,  for a prime
uncolorable formation with contextual curves of different colors.
\bigbreak

Now consider a planar formation $F(T)$ that is minimal, prime and uncolorable. Suppose that it has
contextual curves of the same color. Then it has curve count five by the above reasoning. 
By performing a simple operation on one of the contextual
curves, we obtain a formation $F'$ with contextual curves of different colors. The curve count of $F'$
cannot be four, since four and five have different parity. Therefore the curve count of $F'$ must 
be five or greater and we conclude that $F'$ is factorized.
Similarly, if we begin with 
a formation that is unfactored and incompletable between two curves of different color, then
by operating on one of them we obtain a formation between curves of the same color. The
original curve count is four and the new curve count, being of the same parity, is either less than 
five (and hence solvable) or greater than five (and hence factored). 
This shows that there does not exist a minimal prime uncolorable (incompletable over the 
empty edge) planar trail $F(T).$ If there are uncolorables then there are minimal uncolorables.
Therefore, no minimal uncolorable planar trail is prime. (The trail factors cannot themselves be
uncolorable, since this would contradict mimimality.) But this is a direct contradiction of the
primality principle. Hence the primality principle implies that there are no uncolorable non-empty 
planar trails. 
\bigbreak

\noindent Now consider a minimal uncolorable cubic graph. Such a graph entails
the possible construction of a minimal uncolorable non-empty trail. Drop an edge from the
graph and color the deleted graph. The missing edge cannot have its endpoints on a single curve
(red, blue or alternating) in the corresponding formation since that will allow the filling in of the
missing edge and a coloration of an uncolorable. Therefore we may take the
missing edge to be between two blues. If there is more than one trail factor
between these two blues then we would have a factored minimal uncolorable trail.
Primality implies that there is only one factor. Therefore the Primality
Principle in conjunction with the Parity Lemma implies the non-existence of a
minimal uncolorable cubic graph with a non-empty trail in the coloration of the
deletion (by one edge) of the graph.  The only remaining possibility is that
after deleting one edge, the graph is identical to two curves. The dumbell (See Figure 1) is the
only such graph.   Therefore the primality principle implies the Four Color
Theorem. \vspace{3mm}

\noindent Conversely, assume the Four Color Theorem.  Then indeed there does not
exist a minimal uncolorable non-empty planar prime trail (since it by definition
implies an uncolorable plane cubic graph with no isthmus). Hence the statement
of the Primality Principle is true.  This completes the proof of the Theorem.//
\vspace{3mm}

\noindent {\bf Remark.} This Theorem constitutes a reformulation of  the Four
Color Theorem, in terms of  the Primality Principle.  This reformulation takes
the coloring problem into a new domain. In the work of G. Spencer-Brown  this
reformulation has been investigated in great depth. The capstone of this work is
an algorithm called the {\em parity pass} (\cite{SB} pp. 182-183) that is intended to extend
formations across an uncompleted five-region whenever the given formation does
not already solve by simple operations. Spencer-Brown has stated repeatedly that
this approach gives a proof of the Four-Color Theorem. It is not the purpose of
this paper to give full review of that work. We recommend that the
reader consult Spencer-Brown \cite{SB}. \vspace{3mm}

\section{The Parity Pass}
We shall say that a formaton is {\em $1$-deficient} relative to a graph $G$ if it formates all but
one edge of $G.$
We shall say that a formation is {\em planar uncolorable $1$-deficient} if the underlying graph is
uncolorable if one adds a single designated edge to it. In discussing the Primality Principle in the
last section, we have considered situations where a graph is formated all except for a single edge.
In that section we showed that the Four Color Theorem was equivalent to the  
Primality Principle which states that one cannot build planar uncolorable $1$-deficient
formations by combining colorable trail factors.
Experience in working with the calculus of  formations tends to bolster one's belief in this principle.
There is another  approach to coloring, also due to Spencer-Brown that sheds light on this issue.  It
is well  known since Kempe
\cite{Kempe} that if one could give an algorithm that would color a cubic map in the  plane when a
coloring was given at all but one five-sided region, then any map could be colored with four colors. In
this section, we describe an algorithm (The Parity Pass) that is designed to handle  the five region in
the context of formations. Spencer-Brown asserts that  {\em given a planar formation that is 
$1$-deficient at a five-region, it is either completable by simple operations, or at some stage in the
parity pass algorithm the resulting formation is completable by simple operations}. We refer the reader to \cite{SB}
for more details about the 
context and possible proof of this algorithm.  The purpose of this section is to give a condensed
description of the Parity Pass, and to urge the reader to try it out on ``hard" examples.
\bigbreak

View Figure 17.  This Figure contains a complete diagrammatic summary of the Parity Pass. There is
an initial diagram and five successive transformations $(A,B,C,D,E)$ to related diagrams. The last 
diagram is locally identical to the first diagram.  Each transformation consists in a single 
idemposition of a closed curve on the given formation. In some cases this curve goes through one of the
empty edges, coloring it, while transforming one of the edges at the five-region into an empty
edge. This is a {\em complex operation}. In other cases the transformation is a simple operation
on the given formation. In fact $A$, $C$ and $D$ are complex operations, while $B$ and $E$ are simple
operations. Each operation can be performed if the given formation is not completable by simple 
operations at the five region. Conversely, if one of the operations in the parity pass cannot be
performed, then that starting configuration can be solved by simple operations. We will not prove these
statements here, but we will give a worked example after some further discussion.
\bigbreak

At each stage of Figure 17 we have indicated with small dark circles
the edges along which the idemposition is to be performed to get to the next stage. Specifically,
performing step $A$ requires an idemposition in blue along an alternating curve plus
the drawing of this curve across the missing edge and the cancellation of a blue edge by the operating
curve. The existence of this idemposition is required to perform step $A$. Step $B$ entails idemposition
in red along a purple/blue alternator. This is the same as following the indicated blue curve with a red
idemposition. It is required that the indicated blue curve is distinct from the other blue curve
indicated in the local diagram. Step $C$ entails idemposition in purple along a blue/red alternator.
Step $D$ entails idemposition in blue along a red curve and demands that the two local red segments are
part of one curve. Step $E$ entails purple idemposition along a blue/red alternator that must be
distinct from the other alternator locally indicated. If all steps of the parity pass can be
performed, then one returns to a local configuration at the five region that is the same as the
starting position. 
\bigbreak 

In a given example the reader can deduce from each 
transformed diagram the locus of the putative operation that produces it.
This locus, and the type of operation can be deduced by comparing the changes between the diagram
and its transform. We also leave to the reader the verification that 
{\em when a transform cannot be accomplished, then the domain formation can be completed by simple operations
originating at the five region.}
\bigbreak
 
It is a fascinating exercise to 
perform this algorithm on examples. Success consists in being {\em unable} to apply one of the four
operations of the parity pass, since this inability implicates a solvable formation.
Figures 18 and 19 give two examples for the reader to examine. These exercises involve 
quite a bit of diagrammatic work, but it is worth the effort. In Figures 18.1 to 18.4 we show the work
involved in solving the example of Figure 18 by the parity pass algorithm. In this case one can apply
steps $A$, $B$ and $C$ of the algorithm. Step $D$ cannot be applied and one can see (Figure
18.4) that the formation at this stage can be solved by simple operations. To see this, examime Figure 18.4 and note that
after a simple operaton on the curve $T$, one can place a purple curve around the five region that fills in the missing
edges and does not cancel any remaining edges.
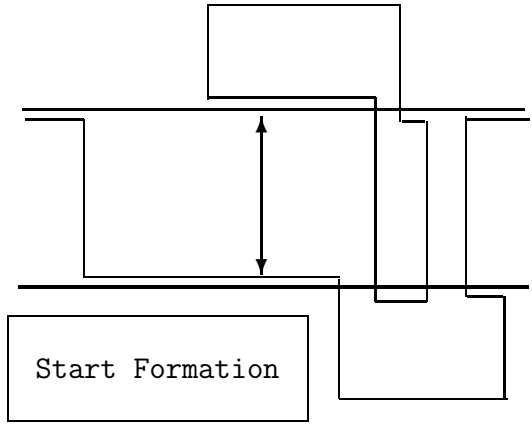
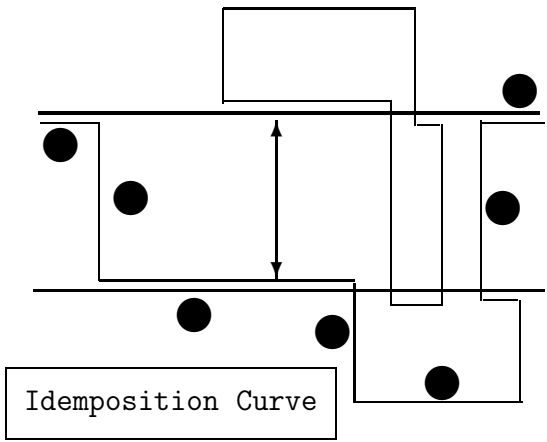
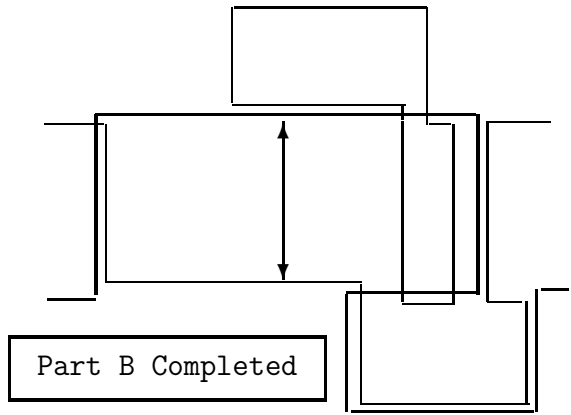
\bigbreak
  
{\tt    \setlength{\unitlength}{0.92pt}
\begin{picture}(422,347)
\thinlines    \put(66,50){\circle*{8}}
              \put(113,62){\circle*{8}}
              \put(123,82){\circle*{8}}
              \put(391,142){\circle*{8}}
              \put(363,154){\circle*{8}}
              \put(389,182){\circle*{8}}
              \put(193,145){\circle*{8}}
              \put(132,143){\circle*{8}}
              \put(101,157){\circle*{8}}
              \put(112,182){\circle*{8}}
              \put(126,200){\circle*{8}}
              \put(395,266){\circle*{8}}
              \put(332,266){\circle*{8}}
              \put(287,279){\circle*{8}}
              \put(269,304){\circle*{8}}
              \put(153,289){\circle*{8}}
              \put(124,277){\circle*{8}}
              \put(85,263){\circle*{8}}
              \put(44,274){\circle*{8}}
              \put(28,289){\circle*{8}}
\thicklines   \put(10,296){\line(1,0){160}}
              \put(12,255){\line(1,0){160}}
\thinlines    \put(90,337){\line(0,-1){37}}
              \put(90,299){\line(1,0){40}}
              \put(130,298){\line(0,-1){41}}
              \put(130,258){\line(1,0){40}}
              \put(50,261){\vector(0,1){29}}
              \put(50,290){\vector(0,-1){29}}
              \put(291,49){\vector(0,-1){29}}
              \put(291,27){\vector(0,1){23}}
              \put(370,270){\vector(0,1){23}}
              \put(370,292){\vector(0,-1){29}}
\thicklines   \put(250,296){\line(1,0){160}}
              \put(251,257){\line(1,0){159}}
\thinlines    \put(251,299){\line(1,0){40}}
              \put(292,299){\line(0,-1){40}}
              \put(292,260){\line(1,0){120}}
              \put(330,337){\line(0,-1){37}}
              \put(330,300){\line(1,0){80}}
\thicklines   \put(178,279){\vector(1,0){62}}
              \put(51,136){\line(1,0){40}}
              \put(91,136){\line(0,1){40}}
              \put(91,176){\line(1,0){120}}
\thinlines    \put(51,176){\line(1,0){37}}
              \put(88,176){\line(0,-1){38}}
              \put(88,138){\line(1,0){123}}
              \put(131,216){\line(0,-1){37}}
              \put(131,179){\line(1,0){80}}
              \put(172,149){\vector(0,1){23}}
              \put(172,171){\vector(0,-1){29}}
              \put(250,176){\line(1,0){160}}
\thicklines   \put(250,136){\line(1,0){160}}
              \put(329,215){\line(0,-1){42}}
              \put(330,173){\line(1,0){80}}
\thinlines    \put(291,147){\vector(0,1){23}}
              \put(291,169){\vector(0,-1){29}}
              \put(370,145){\vector(0,1){23}}
              \put(370,167){\vector(0,-1){29}}
\thicklines   \put(10,157){\vector(1,0){36}}
              \put(212,157){\vector(1,0){36}}
\thinlines    \put(51,56){\line(1,0){119}}
\thicklines   \put(50,17){\line(1,0){160}}
\thinlines    \put(171,56){\line(0,-1){44}}
              \put(171,12){\line(1,0){39}}
\thicklines   \put(130,99){\line(0,-1){38}}
              \put(130,61){\line(1,0){80}}
\thinlines    \put(91,28){\vector(0,1){23}}
              \put(91,50){\vector(0,-1){29}}
\thicklines   \put(249,55){\line(1,0){160}}
              \put(250,15){\line(1,0){160}}
\thinlines    \put(330,97){\line(0,-1){37}}
              \put(330,59){\line(1,0){40}}
              \put(370,59){\line(0,-1){41}}
              \put(370,19){\line(1,0){40}}
\thicklines   \put(211,38){\vector(1,0){36}}
              \put(11,37){\vector(1,0){36}}
\thinlines    \put(198,273){\makebox(25,24){A}}
              \put(16,150){\makebox(23,26){B}}
              \put(221,152){\makebox(22,23){C}}
              \put(19,31){\makebox(22,22){D}}
              \put(218,32){\makebox(21,21){E}}
              \put(15,298){\framebox(17,16){r}}
              \put(14,259){\framebox(17,16){r}}
              \put(257,277){\framebox(17,16){r}}
              \put(256,238){\framebox(17,16){r}}
              \put(56,117){\framebox(17,16){r}}
              \put(259,137){\framebox(17,16){r}}
              \put(333,195){\framebox(17,16){r}}
              \put(134,79){\framebox(17,16){r}}
              \put(57,19){\framebox(17,16){r}}
              \put(255,56){\framebox(17,16){r}}
              \put(256,17){\framebox(17,16){r}}
              \put(92,316){\framebox(16,15){b}}
              \put(333,319){\framebox(16,15){b}}
              \put(294,270){\framebox(16,15){b}}
              \put(58,177){\framebox(16,15){b}}
              \put(133,195){\framebox(16,15){b}}
              \put(260,178){\framebox(16,15){b}}
              \put(59,58){\framebox(16,15){b}}
              \put(333,77){\framebox(16,15){b}}
\end{picture}}

\begin{center}{\bf Figure 17 - The Parity Pass} \end{center}
\bigbreak

  {\tt    \setlength{\unitlength}{0.92pt}
\begin{picture}(362,131)
\thinlines    \put(42,57){\vector(0,1){28}}
              \put(42,85){\vector(0,-1){28}}
\thicklines   \put(1,51){\line(1,0){359}}
              \put(1,92){\line(1,0){360}}
\thinlines    \put(80,128){\line(0,-1){41}}
              \put(221,128){\line(-1,0){141}}
              \put(222,86){\line(0,1){42}}
              \put(200,86){\line(1,0){22}}
              \put(201,46){\line(0,1){41}}
              \put(241,47){\line(-1,0){40}}
              \put(241,85){\line(0,-1){37}}
              \put(283,86){\line(-1,0){42}}
              \put(282,42){\line(0,1){43}}
              \put(322,43){\line(-1,0){40}}
              \put(322,2){\line(0,1){41}}
              \put(162,3){\line(1,0){161}}
              \put(161,45){\line(0,-1){41}}
              \put(121,46){\line(1,0){40}}
              \put(121,87){\line(0,-1){40}}
              \put(81,87){\line(1,0){40}}
\end{picture}}  

\begin{center} {\bf Figure 18 - Culprit Number One} \end{center} \bigbreak

{\tt    \setlength{\unitlength}{0.92pt}
\begin{picture}(238,543)
\thicklines   \put(221,314){\circle*{14}}
              \put(195,314){\circle*{14}}
              \put(171,314){\circle*{14}}
              \put(146,313){\circle*{14}}
              \put(117,273){\circle*{14}}
              \put(84,242){\circle*{14}}
              \put(55,271){\circle*{14}}
              \put(27,314){\circle*{14}}
              \put(44,239){\vector(0,1){63}}
              \put(44,270){\vector(0,-1){33}}
\thinlines    \put(98,309){\line(1,0){29}}
              \put(127,309){\line(0,-1){74}}
              \put(127,236){\line(1,0){26}}
              \put(153,235){\line(0,-1){49}}
              \put(152,186){\line(1,0){69}}
              \put(220,186){\line(0,1){42}}
              \put(219,228){\line(-1,0){14}}
              \put(204,228){\line(0,1){74}}
              \put(204,302){\line(-1,0){14}}
              \put(190,302){\line(0,-1){74}}
              \put(190,228){\line(-1,0){21}}
              \put(169,228){\line(0,1){74}}
              \put(168,302){\line(1,0){11}}
              \put(179,302){\line(0,1){48}}
              \put(178,350){\line(-1,0){79}}
              \put(99,349){\line(0,-1){39}}
\thicklines   \put(20,305){\line(1,0){206}}
              \put(19,232){\line(1,0){206}}
\thinlines    \put(44,126){\line(-1,0){24}}
              \put(201,126){\line(1,0){26}}
\thicklines   \put(13,413){\line(1,0){206}}
              \put(20,57){\line(1,0){206}}
              \put(20,130){\line(1,0){206}}
              \put(14,486){\line(1,0){206}}
\thinlines    \put(177,127){\line(1,0){9}}
              \put(166,134){\line(0,-1){6}}
              \put(95,134){\line(1,0){71}}
              \put(44,61){\line(0,1){65}}
              \put(123,61){\line(-1,0){79}}
\thicklines   \put(117,64){\vector(0,1){63}}
              \put(117,95){\vector(0,-1){33}}
\thinlines    \put(123,61){\line(1,0){26}}
              \put(149,60){\line(0,-1){49}}
              \put(148,12){\line(1,0){69}}
              \put(216,11){\line(0,1){42}}
              \put(215,53){\line(-1,0){14}}
              \put(201,53){\line(0,1){74}}
              \put(187,127){\line(0,-1){74}}
              \put(187,53){\line(-1,0){21}}
              \put(166,53){\line(0,1){74}}
              \put(176,126){\line(0,1){48}}
              \put(175,174){\line(-1,0){79}}
              \put(95,173){\line(0,-1){39}}
              \put(92,529){\line(0,-1){39}}
              \put(171,531){\line(-1,0){79}}
              \put(173,483){\line(0,1){48}}
              \put(162,483){\line(1,0){11}}
              \put(163,409){\line(0,1){74}}
              \put(185,409){\line(-1,0){21}}
              \put(184,483){\line(0,-1){74}}
              \put(198,482){\line(-1,0){14}}
              \put(198,409){\line(0,1){74}}
              \put(212,409){\line(-1,0){14}}
              \put(213,367){\line(0,1){42}}
              \put(145,367){\line(1,0){69}}
              \put(146,416){\line(0,-1){49}}
              \put(120,417){\line(1,0){26}}
              \put(120,490){\line(0,-1){74}}
              \put(91,490){\line(1,0){29}}
\thicklines   \put(38,451){\vector(0,-1){33}}
              \put(38,420){\vector(0,1){63}}
\thinlines    \put(14,362){\framebox(123,43){Start Formation}}
              \put(10,189){\framebox(135,29){Idemposition Curve}}
              \put(16,17){\framebox(127,25){Part A  Completed}}
\end{picture}}

\begin{center} {\bf Figure 18.1 - Stage A of Parity Pass Applied to Culprit Number One.} \end{center}
\bigbreak

{\tt    \setlength{\unitlength}{0.92pt}
\begin{picture}(257,561)
\thicklines   \put(224,328){\circle*{14}}
              \put(217,280){\circle*{14}}
              \put(192,208){\circle*{14}}
              \put(147,229){\circle*{14}}
              \put(90,236){\circle*{14}}
              \put(64,284){\circle*{14}}
              \put(35,306){\circle*{14}}
\thinlines    \put(171,324){\line(0,-1){84}}
              \put(102,324){\line(1,0){69}}
\thicklines   \put(24,246){\line(1,0){210}}
\thinlines    \put(51,315){\line(-1,0){24}}
              \put(208,315){\line(1,0){26}}
\thicklines   \put(26,319){\line(1,0){206}}
\thinlines    \put(182,314){\line(1,0){9}}
              \put(51,250){\line(0,1){65}}
              \put(130,250){\line(-1,0){79}}
\thicklines   \put(124,253){\vector(0,1){63}}
              \put(124,284){\vector(0,-1){33}}
\thinlines    \put(130,250){\line(1,0){26}}
              \put(156,249){\line(0,-1){49}}
              \put(156,200){\line(1,0){69}}
              \put(224,200){\line(0,1){42}}
              \put(223,242){\line(-1,0){14}}
              \put(208,242){\line(0,1){74}}
              \put(192,314){\line(0,-1){74}}
              \put(192,240){\line(-1,0){21}}
              \put(181,314){\line(0,1){48}}
              \put(181,362){\line(-1,0){79}}
              \put(102,362){\line(0,-1){39}}
              \put(13,185){\framebox(135,29){Idemposition Curve}}
              \put(10,378){\framebox(123,43){Start Formation}}
              \put(92,549){\line(0,-1){39}}
              \put(171,549){\line(-1,0){79}}
              \put(171,501){\line(0,1){48}}
              \put(182,427){\line(-1,0){21}}
              \put(182,501){\line(0,-1){74}}
              \put(198,429){\line(0,1){74}}
              \put(213,429){\line(-1,0){14}}
              \put(214,387){\line(0,1){42}}
              \put(146,387){\line(1,0){69}}
              \put(146,436){\line(0,-1){49}}
              \put(120,437){\line(1,0){26}}
\thicklines   \put(114,471){\vector(0,-1){33}}
              \put(114,440){\vector(0,1){63}}
\thinlines    \put(120,437){\line(-1,0){79}}
              \put(41,437){\line(0,1){65}}
              \put(172,501){\line(1,0){9}}
\thicklines   \put(16,506){\line(1,0){206}}
\thinlines    \put(198,502){\line(1,0){26}}
              \put(41,502){\line(-1,0){24}}
\thicklines   \put(14,433){\line(1,0){210}}
\thinlines    \put(56,130){\line(-1,0){24}}
              \put(214,131){\line(1,0){26}}
              \put(190,130){\line(1,0){9}}
              \put(179,138){\line(0,-1){6}}
              \put(109,138){\line(1,0){71}}
              \put(57,65){\line(0,1){65}}
              \put(136,65){\line(-1,0){79}}
\thicklines   \put(130,68){\vector(0,1){63}}
              \put(130,99){\vector(0,-1){33}}
\thinlines    \put(136,65){\line(1,0){26}}
              \put(162,64){\line(0,-1){49}}
              \put(162,15){\line(1,0){69}}
              \put(230,15){\line(0,1){42}}
              \put(228,57){\line(-1,0){14}}
              \put(214,57){\line(0,1){74}}
              \put(200,130){\line(0,-1){74}}
              \put(200,56){\line(-1,0){21}}
              \put(179,57){\line(0,1){74}}
              \put(189,130){\line(0,1){48}}
              \put(189,178){\line(-1,0){79}}
              \put(109,178){\line(0,-1){39}}
\thicklines   \put(53,134){\line(0,-1){74}}
              \put(53,58){\line(-1,0){20}}
              \put(53,134){\line(1,0){156}}
              \put(209,61){\line(-1,0){53}}
              \put(156,60){\line(0,-1){48}}
              \put(158,12){\line(1,0){75}}
              \put(234,12){\line(0,1){50}}
              \put(236,62){\line(1,0){11}}
              \put(210,134){\line(0,-1){74}}
              \put(18,17){\framebox(129,25){Part B Completed}}
\thinlines    \put(92,511){\line(1,0){69}}
              \put(161,511){\line(0,-1){84}}
\end{picture}}

\begin{center} {\bf Figure 18.2 - Stage B of Parity Pass Applied to Culprit Number One.} \end{center}
\bigbreak

{\tt    \setlength{\unitlength}{0.92pt}
\begin{picture}(259,504)
\thicklines   \put(10,17){\framebox(131,23){Part C Completed}}
              \put(136,71){\vector(0,1){25}}
              \put(136,83){\vector(0,-1){32}}
              \put(71,85){\vector(0,-1){32}}
              \put(71,73){\vector(0,1){25}}
              \put(185,131){\line(-1,0){61}}
              \put(185,106){\line(0,1){25}}
              \put(215,105){\line(-1,0){31}}
              \put(215,44){\line(0,1){61}}
              \put(177,44){\line(1,0){37}}
              \put(176,100){\line(0,-1){56}}
              \put(124,100){\line(1,0){52}}
              \put(123,130){\line(0,-1){30}}
              \put(232,49){\line(1,0){17}}
              \put(230,12){\line(0,1){37}}
              \put(146,12){\line(1,0){84}}
              \put(145,47){\line(0,-1){35}}
              \put(146,47){\line(-1,0){104}}
\thinlines    \put(149,47){\line(0,-1){31}}
              \put(212,102){\line(1,0){29}}
              \put(211,47){\line(0,1){55}}
              \put(227,47){\line(-1,0){15}}
              \put(227,17){\line(0,1){29}}
              \put(148,16){\line(1,0){79}}
              \put(196,47){\line(-1,0){47}}
              \put(197,103){\line(0,-1){56}}
              \put(44,103){\line(1,0){153}}
              \put(78,274){\circle*{12}}
              \put(110,294){\circle*{12}}
              \put(177,276){\circle*{12}}
              \put(165,229){\circle*{12}}
              \put(166,185){\circle*{12}}
              \put(105,190){\circle*{12}}
              \put(59,229){\circle*{12}}
\thicklines   \put(204,267){\line(0,-1){74}}
              \put(230,195){\line(1,0){11}}
              \put(228,145){\line(0,1){50}}
              \put(152,145){\line(1,0){75}}
              \put(150,193){\line(0,-1){48}}
              \put(203,194){\line(-1,0){53}}
              \put(47,267){\line(1,0){156}}
              \put(47,191){\line(-1,0){20}}
              \put(47,267){\line(0,-1){74}}
\thinlines    \put(102,311){\line(0,-1){39}}
              \put(182,311){\line(-1,0){79}}
              \put(183,263){\line(0,1){48}}
              \put(173,190){\line(0,1){74}}
              \put(194,189){\line(-1,0){21}}
              \put(194,263){\line(0,-1){74}}
              \put(208,190){\line(0,1){74}}
              \put(223,190){\line(-1,0){14}}
              \put(224,148){\line(0,1){42}}
              \put(156,148){\line(1,0){69}}
              \put(156,197){\line(0,-1){49}}
              \put(130,198){\line(1,0){26}}
\thicklines   \put(125,236){\vector(0,-1){33}}
\thinlines    \put(130,198){\line(-1,0){79}}
              \put(51,198){\line(0,1){65}}
              \put(102,271){\line(1,0){71}}
              \put(173,271){\line(0,-1){6}}
              \put(184,263){\line(1,0){9}}
              \put(208,264){\line(1,0){26}}
              \put(50,263){\line(-1,0){24}}
\thicklines   \put(125,219){\vector(0,1){43}}
              \put(121,400){\vector(0,1){43}}
\thinlines    \put(11,144){\framebox(135,29){Idemposition Curve}}
              \put(17,320){\framebox(123,43){Start Formation}}
              \put(46,444){\line(-1,0){24}}
              \put(204,445){\line(1,0){26}}
              \put(180,444){\line(1,0){9}}
              \put(169,452){\line(0,-1){6}}
              \put(98,452){\line(1,0){71}}
              \put(47,379){\line(0,1){65}}
              \put(126,379){\line(-1,0){79}}
\thicklines   \put(121,417){\vector(0,-1){33}}
\thinlines    \put(126,379){\line(1,0){26}}
              \put(152,378){\line(0,-1){49}}
              \put(152,329){\line(1,0){69}}
              \put(220,329){\line(0,1){42}}
              \put(219,371){\line(-1,0){14}}
              \put(204,371){\line(0,1){74}}
              \put(190,444){\line(0,-1){74}}
              \put(190,370){\line(-1,0){21}}
              \put(169,371){\line(0,1){74}}
              \put(179,444){\line(0,1){48}}
              \put(178,492){\line(-1,0){79}}
              \put(98,492){\line(0,-1){39}}
\thicklines   \put(43,448){\line(0,-1){74}}
              \put(43,372){\line(-1,0){20}}
              \put(43,448){\line(1,0){156}}
              \put(199,375){\line(-1,0){53}}
              \put(146,374){\line(0,-1){48}}
              \put(148,326){\line(1,0){75}}
              \put(224,326){\line(0,1){50}}
              \put(226,376){\line(1,0){11}}
              \put(200,448){\line(0,-1){74}}
\end{picture}}

\begin{center} {\bf Figure 18.3 - Stage C of Parity Pass Applied to Culprit Number One.} \end{center}
\bigbreak

{\tt    \setlength{\unitlength}{0.92pt}
\begin{picture}(277,398)
\thicklines   \put(10,10){\framebox(248,34){solution at a-b-c-d-e.}}
\thinlines    \put(79,173){\framebox(11,13){e}}
\thicklines   \put(154,326){\vector(0,1){25}}
              \put(154,338){\vector(0,-1){32}}
              \put(89,340){\vector(0,-1){32}}
              \put(89,328){\vector(0,1){25}}
              \put(203,386){\line(-1,0){61}}
              \put(203,361){\line(0,1){25}}
              \put(233,360){\line(-1,0){31}}
              \put(233,299){\line(0,1){61}}
              \put(195,299){\line(1,0){37}}
              \put(194,355){\line(0,-1){56}}
              \put(142,355){\line(1,0){52}}
              \put(141,385){\line(0,-1){30}}
              \put(250,304){\line(1,0){17}}
              \put(248,267){\line(0,1){37}}
              \put(164,267){\line(1,0){84}}
              \put(163,302){\line(0,-1){35}}
              \put(164,302){\line(-1,0){104}}
\thinlines    \put(167,302){\line(0,-1){31}}
              \put(230,357){\line(1,0){29}}
              \put(229,302){\line(0,1){55}}
              \put(245,302){\line(-1,0){15}}
              \put(245,272){\line(0,1){29}}
              \put(166,271){\line(1,0){79}}
              \put(214,302){\line(-1,0){47}}
              \put(215,358){\line(0,-1){56}}
              \put(62,358){\line(1,0){153}}
\thicklines   \put(14,268){\framebox(143,23){Start Formation}}
\thinlines    \put(24,227){\line(1,0){153}}
              \put(177,227){\line(0,-1){56}}
              \put(176,171){\line(-1,0){47}}
              \put(128,140){\line(1,0){79}}
              \put(207,141){\line(0,1){29}}
              \put(207,171){\line(-1,0){15}}
              \put(191,171){\line(0,1){55}}
              \put(192,226){\line(1,0){29}}
              \put(129,171){\line(0,-1){31}}
\thicklines   \put(125,171){\line(-1,0){104}}
              \put(125,171){\line(0,-1){35}}
              \put(210,136){\line(0,1){37}}
              \put(212,173){\line(1,0){17}}
              \put(103,254){\line(0,-1){30}}
              \put(104,224){\line(1,0){52}}
              \put(156,224){\line(0,-1){56}}
              \put(157,168){\line(1,0){37}}
              \put(195,168){\line(0,1){61}}
              \put(195,229){\line(-1,0){31}}
              \put(165,230){\line(0,1){25}}
              \put(165,255){\line(-1,0){61}}
              \put(51,197){\vector(0,1){25}}
              \put(51,209){\vector(0,-1){32}}
              \put(116,207){\vector(0,-1){32}}
              \put(116,195){\vector(0,1){25}}
              \put(134,218){\circle*{6}}
              \put(122,197){\circle*{6}}
              \put(123,176){\circle*{6}}
              \put(12,80){\framebox(242,43){Idemposition Curve does not exist.}}
              \put(125,135){\line(1,0){85}}
              \put(168,238){\framebox(15,16){T}}
              \put(38,193){\framebox(12,13){a}}
              \put(78,211){\framebox(12,13){b}}
              \put(106,228){\framebox(11,13){c}}
              \put(101,190){\framebox(13,14){d}}
              \put(10,42){\framebox(245,35){Simple operation on curve T allows}}
\end{picture}}

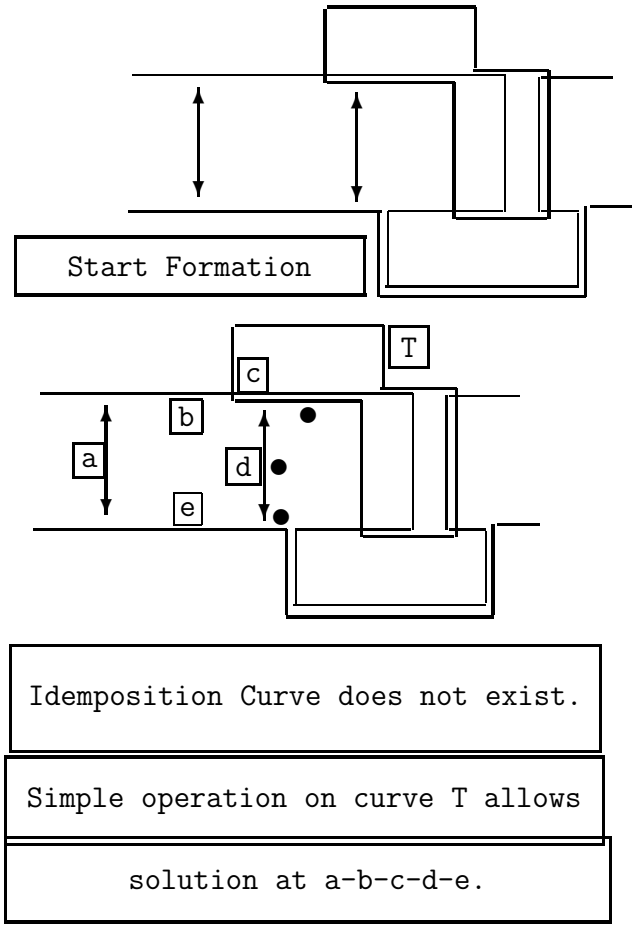
\begin{center} {\bf Figure 18.4 - Stage D of Parity Pass Applied to Culprit Number One.
Idemposition cannot be performed and Formation is Completable at the Five Region.} \end{center}
\bigbreak

{\tt    \setlength{\unitlength}{0.92pt}
\begin{picture}(456,303)
\thicklines   \put(1,99){\line(1,0){452}}
              \put(3,219){\line(1,0){452}}
\thinlines    \put(207,102){\line(1,0){21}}
              \put(231,222){\line(-1,-5){24}}
              \put(204,222){\line(1,0){27}}
              \put(204,277){\line(0,-1){55}}
\thicklines   \put(21,199){\vector(0,-1){94}}
              \put(21,106){\vector(0,1){109}}
\thinlines    \put(241,214){\line(1,0){24}}
              \put(264,213){\line(-1,-3){37}}
              \put(382,222){\line(-1,0){16}}
              \put(382,277){\line(0,-1){55}}
              \put(80,215){\line(1,0){26}}
              \put(80,240){\line(0,-1){26}}
              \put(185,104){\line(0,1){118}}
              \put(108,17){\line(1,0){157}}
              \put(97,3){\line(0,1){91}}
              \put(319,3){\line(-1,0){222}}
              \put(319,180){\line(0,1){35}}
              \put(63,105){\line(0,1){109}}
              \put(382,277){\line(-1,0){179}}
              \put(366,259){\line(0,-1){37}}
              \put(242,259){\line(1,0){124}}
              \put(241,215){\line(0,1){44}}
              \put(42,214){\line(1,0){21}}
              \put(42,300){\line(0,-1){86}}
              \put(418,300){\line(-1,0){375}}
              \put(419,215){\line(0,1){85}}
              \put(439,214){\line(-1,0){19}}
              \put(440,181){\line(0,1){33}}
              \put(319,180){\line(1,0){121}}
              \put(298,215){\line(1,0){21}}
              \put(297,93){\line(0,1){122}}
              \put(318,93){\line(-1,0){20}}
              \put(319,3){\line(0,1){90}}
              \put(107,94){\line(-1,0){9}}
              \put(108,17){\line(0,1){77}}
              \put(265,94){\line(0,-1){76}}
              \put(245,94){\line(1,0){20}}
              \put(245,61){\line(0,1){32}}
              \put(142,60){\line(1,0){103}}
              \put(142,103){\line(0,-1){43}}
              \put(125,103){\line(1,0){17}}
              \put(125,215){\line(0,-1){112}}
              \put(141,216){\line(-1,0){16}}
              \put(163,105){\line(-1,5){22}}
              \put(185,104){\line(-1,0){22}}
              \put(163,222){\line(1,0){22}}
              \put(163,240){\line(0,-1){18}}
              \put(80,240){\line(1,0){82}}
              \put(84,104){\line(1,5){22}}
              \put(63,104){\line(1,0){21}}
\end{picture}}

\begin{center} {\bf Figure 19 - Culprit Number Two} \end{center}
\bigbreak

In these examples, Culprit Number One (Figure 18) solves via parity pass after application of $A$,
$B$ and $C.$  Culprit Number Two (Figure 19) will solve via parity pass after application of 
$A$ and $B.$  Culprit Number Two is an example of a ``good try" at making a factorized minimal
uncolorable in the plane.  As we mentioned, Spencer-Brown asserts that either the 
parity pass solves any five region extension problem, or the problem could have been solved by simple 
operations at the outset. 

\section {The Penrose Formula}
Roger Penrose \cite{P} gives a formula for computing the number of proper 
edge 3-colorings of a plane cubic graph $G.$
In this formula each vertex is associated with the ``epsilon" tensor 
$$P_{ijk}=\sqrt{-1}\epsilon_{ijk}$$ 
\noindent as shown in Figure 20.

{\tt    \setlength{\unitlength}{0.92pt}
\begin{picture}(227,222)
\thinlines    \put(212,81){\makebox(14,18){r}}
              \put(131,79){\makebox(16,19){b}}
              \put(70,75){\makebox(16,19){b}}
              \put(15,76){\makebox(14,18){r}}
              \put(163,39){\line(0,-1){38}}
              \put(128,82){\line(5,-6){35}}
\thicklines   \put(168,41){\line(0,-1){39}}
              \put(209,82){\line(-1,-1){41}}
\thinlines    \put(54,43){\line(0,-1){40}}
              \put(88,78){\line(-1,-1){34}}
\thicklines   \put(48,40){\line(0,-1){37}}
              \put(8,82){\line(1,-1){40}}
\thinlines    \put(83,202){\makebox(16,19){b}}
              \put(39,98){\makebox(17,20){p}}
              \put(1,203){\makebox(14,18){r}}
              \put(123,202){\makebox(16,19){b}}
              \put(159,98){\makebox(17,20){p}}
              \put(202,203){\makebox(14,18){r}}
              \put(183,148){\makebox(20,19){$~-\sqrt{-1}$}}
\thicklines   \put(129,201){\line(1,-1){39}}
              \put(167,161){\line(1,1){40}}
              \put(167,161){\line(0,-1){41}}
              \put(168,161){\circle*{18}}
              \put(67,150){\makebox(21,21){$~+\sqrt{-1}$}}
              \put(48,161){\circle*{18}}
              \put(47,161){\line(0,-1){41}}
              \put(47,161){\line(1,1){40}}
              \put(9,201){\line(1,-1){39}}
\end{picture}}

\begin{center} {\bf Figure 20 - Epsilon Tensor} \end{center}
\bigbreak

One takes the colors from the set $\{1,2,3\}$ and the tensor $\epsilon_{ijk}$ takes value $1$
for $ijk = 123, 231,312$ and $-1$ for $ijk = 132, 321, 213.$  The tensor is $0$ when $ijk$ is not 
a permutation of $123.$ One then evaluates the graph $G$ by taking the sum over all possible color
assigments to its edges of the products of the $P_{ijk}$ associated with its nodes. Call this 
evaluation $[G].$ 
\bigbreak

\noindent {\bf Theorem (Penrose).} If $G$ is a planar cubic graph, then $[G]$, as defined above, is equal to the number of
distinct proper colorings of the edges of $G$ with three colors (so that every vertex sees three colors at its edges).
\bigbreak

\noindent {\bf Proof.} It follows from the above description that only proper colorings of $G$ contribute to 
the summation $[G],$ and that each such coloring contributes a product of $\pm \sqrt{-1}$ from the tensor
evaluations at the nodes of the graph. In order to see that $[G]$ is equal to the number of 
colorings for a plane graph, one must see that each such contribution is equal to $+1.$
The proof of this assertion is given in Figure 21 where we see that in a formation for a coloring 
each bounce contributes $+1 = -\sqrt{-1}\sqrt{-1}$ while each crossing contributes $-1.$ Since there
are an even number of crossings among the curves in the formation, it follows that the total 
product is equal to $+1$. This completes the proof of the Penrose Theorem.
\bigbreak

{\tt    \setlength{\unitlength}{0.92pt}
\begin{picture}(177,286)
\thinlines    \put(153,41){\framebox(22,20){$-1$}}
              \put(154,161){\framebox(22,20){$-1$}}
              \put(151,240){\framebox(24,23){$+1$}}
              \put(68,66){\makebox(15,17){b}}
              \put(4,43){\makebox(14,16){r}}
              \put(49,185){\makebox(15,17){b}}
              \put(5,162){\makebox(14,16){r}}
              \put(89,44){\makebox(18,18){$~-\sqrt{-1}$}}
              \put(48,19){\makebox(18,18){$~-\sqrt{-1}$}}
              \put(92,140){\makebox(15,17){$~\sqrt{-1}$}}
              \put(32,139){\makebox(15,17){$~\sqrt{-1}$}}
              \put(44,46){\line(0,-1){45}}
              \put(84,46){\line(-1,0){39}}
              \put(84,82){\line(0,-1){36}}
\thicklines   \put(4,42){\line(1,0){120}}
\thinlines    \put(86,165){\line(0,-1){45}}
              \put(46,165){\line(1,0){40}}
              \put(45,201){\line(0,-1){36}}
\thicklines   \put(5,161){\line(1,0){120}}
\thinlines    \put(48,268){\makebox(15,17){b}}
              \put(1,243){\makebox(14,16){r}}
              \put(84,220){\makebox(18,18){$~-\sqrt{-1}$}}
              \put(31,222){\makebox(15,17){$~\sqrt{-1}$}}
              \put(86,248){\line(0,1){35}}
              \put(44,247){\line(1,0){42}}
              \put(44,284){\line(0,-1){36}}
\thicklines   \put(4,242){\line(1,0){120}}
\end{picture}}

\begin{center} {\bf Figure 21 - Cross and Bounce} \end{center}
\bigbreak

It is easy to see from the properties of the epsilon tensor that $[G]$ 
satisfies the recursive identity shown in Figure 22. Here we have that $[O] = 3,$ where $O$ denotes
an isolated curve, and the recursion formula includes graphs with extra crossings as shown in the
Figure. This use of formations gives a vivid access to the theory of the Penrose 
formula.
\bigbreak

{\tt    \setlength{\unitlength}{0.92pt}
\begin{picture}(326,127)
\thicklines   \put(209,55){\makebox(24,22){$-$}}
              \put(84,55){\makebox(24,23){$=$}}
              \put(183,25){\line(1,-1){20}}
              \put(182,106){\line(0,-1){80}}
              \put(201,124){\line(-1,-1){19}}
              \put(139,23){\line(-1,-1){17}}
              \put(140,106){\line(0,-1){81}}
              \put(122,124){\line(1,-1){18}}
              \put(242,5){\line(2,3){80}}
              \put(242,124){\line(2,-3){81}}
              \put(43,42){\circle*{20}}
              \put(42,84){\circle*{20}}
              \put(41,44){\line(1,-1){40}}
              \put(42,44){\line(-1,-1){41}}
              \put(40,84){\line(0,-1){40}}
              \put(41,83){\line(1,1){40}}
              \put(1,124){\line(1,-1){40}}
\end{picture}}

\begin{center} {\bf Figure 22 - Penrose Formula} \end{center}
\bigbreak

\section{The Eliahou-Kryuchkov Conjecture}
The Eliahou-Kryuchkov Conjecture \cite{Kryuchkov, Eliahou} is about ``reassociating" signed trees.
The term {\em reassociation} comes from the algebraic transform of a product $(ab)c$ to a product $a(bc).$
In a non-associative algebra these two terms can represent distinct algebraic elements.
In a tree, a trivalent vertex can be regarded as a representative for an algebraic product in the sense that 
two edge labels
are multiplied to give the third edge label at that vertex. See Figure 23 for an illustration of this pattern. In this
Figure we show how two distinct trees correspond to the two associated products $(ab)c$ and $a(bc).$
\bigbreak
  
{\tt    \setlength{\unitlength}{0.92pt}
\begin{picture}(390,416)
\thicklines   \put(276,185){\framebox(47,17){(ab)c}}
              \put(17,386){\line(1,-1){59}}
              \put(74,329){\line(3,4){45}}
              \put(75,327){\line(0,-1){67}}
              \put(75,327){\circle*{16}}
              \put(10,385){\framebox(19,19){a}}
              \put(111,388){\framebox(20,18){b}}
              \put(60,238){\framebox(35,20){ab}}
              \put(264,388){\framebox(20,18){b}}
              \put(162,385){\framebox(19,19){a}}
              \put(228,329){\line(3,4){45}}
              \put(171,386){\line(1,-1){59}}
              \put(231,325){\line(3,-2){67}}
              \put(299,278){\line(2,3){68}}
              \put(298,279){\line(0,-1){74}}
              \put(359,380){\framebox(21,16){c}}
              \put(300,278){\circle*{16}}
              \put(230,329){\circle*{16}}
              \put(37,182){\line(4,-3){122}}
              \put(159,90){\line(1,1){89}}
              \put(153,180){\line(6,-5){54}}
              \put(158,89){\line(0,-1){56}}
              \put(206,137){\circle*{16}}
              \put(158,90){\circle*{16}}
              \put(29,181){\framebox(19,19){a}}
              \put(143,180){\framebox(20,18){b}}
              \put(239,180){\framebox(21,16){c}}
              \put(137,10){\framebox(45,20){a(bc)}}
\end{picture}}

\begin{center} {\bf Figure 23 - Multiplication, Trees and Associated Products} \end{center}
\bigbreak

The basic reassociation pattern in binary tree form (each vertex is incident to three edges) is  shown
in Figures 23 and 24 . Note that in the lower half of Figure 24 the trees are labeled with the colors $p$, $r$
and
$b$ with the product of any two of these colors equal to the third color. In this case we see that  the
tree diagram has illustrated the identity $(rb)r = pr = b = rp = r(br).$  Thus in this case the
multiplication is associative. If we decide that $rr = pp =bb = 0$ with $0r = r0 = 0b =b0 = 0p = p0 =
0$ then the system
$\{r,b,p,0\}$ is not associative and two  colored trees changed by a reassociation move as indicated in
the Figure may have different coloring properties.  Note that if we have a colored tree (three distinct
colors at a vertex) we can take the two cyclic color orders ($rbp$ clockwise or 
$rpb$ clockwise) as denoting two possible signs ( plus and minus respectively) that can be assigned to the vertices of the tree.
In the context of the conjecture we are about to discuss one considers arbitrary assignments of signs to the vertices of a tree.
The relation with coloring is left out of the game momentarily. 
\smallbreak

Here is a remarkable game! We shall give signs to the vertices of a tree. We allow the reassociation move inside a larger tree only when
the two adjacent vertices in the reassociation are assigned the same sign, and then 
both vertices receive the opposite of this sign after the reassociation. Such moves are called {\it signed reassociation moves.}
\bigbreak 

\noindent {\bf Eliahou-Kryuchkov Conjecture.}  Given any two connected trees (with cubic vertices and, at the ends,
vertices incident to single edges) and the same number of twigs (a {\em twig} is an edge incident to an
end vertex in the tree) , {\it then there exist assignments of signs to the vertices of the two trees so that one signed
tree can be transformed to the other signed tree by a series of signed reassociation moves.} We shall
refer to this conjecture as the EK conjecture.
\bigbreak

It was known to the authors of this conjecture that the four color theorem follows from it.
In \cite{GP} it has been shown that in fact the EK conjecture is equivalent to the four color 
theorem. We mention the EK conjecture here to point out that it implies that {\it one can edge color the two 
trees so that one tree can be obtained from the other by reassociation moves on the colorings as 
shown in Figure 24 using formations.} These reassociation moves on the colorings are particulary nice in that they do not involve changing the colors
only reconfiguring the graph.  The proof of this statement follows directly from the local coloring depicted in Figure 24. There is a
particularly nice pathway of colorings leading from one colored tree to the other. 
\bigbreak

In this way the formations make the nature of the reassociation move 
clear and show how coloring is related to the EK conjecture. 
It is remarkable that the four color theorem is equivalent to this very specific
statement  about coloring trees.
\bigbreak 

We can also see just how the EK conjecture is related to the 
vector cross product reformulation of the four color theorem \cite{VCP}. In the vector cross product
reformulation of the four color theorem, we are given two associated products of the same ordered 
sets of variables. The {\em vector product conjecture} then states that there exist assignments
to the variables from the set of generators $\{i,j,k\}$ of the vector cross product algebra in 
three dimensional space such that each of the two given products is non-zero in this algebra. This is
sufficient to make the two products equal since if they are non-zero then all partial products are
non-zero and hence each product may be viewed in the quaternions. Since the quaternions are
associative, it  follows that the two products are equal. In this sense, the vector product conjecture is actually a conjecture about the structure of
the quaternions. 
\bigbreak

The relationship of the vector product conjecture with graph coloring is
obtained by forming a plane graph consisting of the two trees, tied at their single roots and tied at
their branches by non-intersecting arcs so that the left-most branch of the left tree has the same
variable as the right-most branch of the right tree and the product order in the left tree is
left-to-right, while the product order in the right tree is right-to-left. {\it Solving the equality of the
two products is equivalent to coloring the graph consisting in two tied trees.}
\bigbreak 

Let the two associations of the product of $n$ variables be denoted $L$ and $R.$ 
\bigbreak

\noindent {\bf Proposition.} The EK conjecture implies that there exists a
solution to the equation $L=R$ in the vector cross product algebra plus a series of 
algebraic reassociations taking $L$ to $R$ such that all of the intermediate 
terms in the sequence of reassociations are  non-zero when evaluated as vector cross products. 
\bigbreak

\noindent {\bf Proof.} The proof of this assertion is easy to see
using the formalism of formations by translating signs to colors as we have illustrated in Figure 24, and using the interpretation of products via
trees as shown in Figures 23 and 24. The signs at the vertices are derived from the fact that in the cross product algebra we have $ij = + k$ and 
$ji = - k.$ One replaces $r,b,p$ by $i,j,k.$ Local signs in the partial products in the trees can then be used to decorate the vertices of the tree.
This completes the sketch of the proof.
\bigbreak
 
This extra texture in the vector cross product formulation, and its relationship 
with the quaternions may provide new algebraic insight into the nature of the four color 
theorem.
\bigbreak

  {\tt    \setlength{\unitlength}{0.92pt}
\begin{picture}(406,283)
\thinlines    \put(56,107){\makebox(15,16){b}}
              \put(86,4){\makebox(15,16){b}}
              \put(329,1){\makebox(15,16){b}}
              \put(307,107){\makebox(15,16){b}}
              \put(375,105){\makebox(16,16){r}}
              \put(260,104){\makebox(16,16){r}}
              \put(132,105){\makebox(16,16){r}}
              \put(17,107){\makebox(16,16){r}}
\thicklines   \put(174,64){\vector(1,0){41}}
              \put(215,64){\vector(-1,0){41}}
              \put(209,206){\vector(-1,0){41}}
              \put(169,206){\vector(1,0){40}}
\thinlines    \put(329,220){\makebox(21,20){$-$}}
              \put(20,219){\makebox(20,21){$+$}}
              \put(294,183){\makebox(21,20){$-$}}
              \put(91,181){\makebox(20,21){$+$}}
              \put(324,44){\line(0,-1){41}}
              \put(361,82){\line(-1,-1){37}}
              \put(323,121){\line(1,-1){38}}
\thicklines   \put(324,40){\line(1,1){80}}
              \put(242,121){\line(1,-1){80}}
\thinlines    \put(82,50){\line(0,-1){45}}
              \put(47,85){\line(1,-1){34}}
              \put(82,120){\line(-1,-1){34}}
\thicklines   \put(83,41){\line(1,1){79}}
              \put(3,121){\line(1,-1){80}}
              \put(283,201){\line(0,-1){40}}
              \put(283,281){\line(1,-1){41}}
              \put(283,199){\line(1,1){80}}
              \put(203,281){\line(1,-1){80}}
              \put(82,199){\line(0,-1){39}}
              \put(82,281){\line(-1,-1){39}}
              \put(82,200){\line(1,1){81}}
              \put(1,281){\line(1,-1){81}}
\end{picture}} 
   
\begin{center} {\bf Figure 24 - Signed Reassociation} \end{center}
\bigbreak

\end{document}